\documentclass[twoside]{article}
\usepackage{graphicx,amssymb,mathrsfs,amsmath,color,fancyhdr,amsthm}
\usepackage[all]{xy}
\input cyracc.def

\renewcommand{\baselinestretch}{1.06} 
\catcode`@=11 \long\def\@makefntext#1{\noindent #1}
\newskip\tabcentering \tabcentering=1000pt plus 1000pt minus 1000pt
\def\MCH#1#2{\setbox0=\hbox{\raise#1\hbox{#2}}\smash{\box0}}

\def\dl{\displaystyle}
\let\@oddfoot\@empty  \let\@evenfoot\@empty

\def\@evenhead{}\def\@oddhead{}
\def\@evenhead{\vbox{\hbox to \textwidth{\footnotesize\rm\hbox to
1.0cm{\thepage\hfill} \hfill\hspace{2mm}\footnotesize{
\emph{Xiao J. and Zhao M.}}}}}
\def\@oddhead{\vbox{\hbox to \textwidth{\footnotesize
{\it BGP-Reflection Functors and Lusztig's Symmetries} \hfill{\ } \hfill\hbox to
1cm{\hfill\thepage}}}}

\newenvironment{prof}[1][Proof]{\noindent\textit{#1}\quad }
{\hfill $\Box$\vspace{0.7mm}}

\begin{document}

\newtheorem{theorem}{Theorem}[section]
\newtheorem{lemma}[theorem]{Lemma}
\newtheorem{corollary}[theorem]{Corollary}
\newtheorem{definition}[theorem]{Definition}
\newtheorem{proposition}[theorem]{Proposition}
\newtheorem{conjecture}[theorem]{Conjecture}
 \newtheorem{algorithm}[theorem]{Algorithm}
 \newtheorem{algo}{Algorithm}
  \newtheorem{thm}{Theorem}
  \newtheorem{lem}[thm]{Lemma}
  \newtheorem{prop}[thm]{Proposition}
  \newtheorem{cor}[thm]{Corollary}
  \newtheorem{defn}{Definition}
  \newtheorem{conj}{Conjecture}
\theoremstyle{definition}
\newtheorem{remark}[theorem]{Remark} 
\newtheorem{problem}[theorem]{Problem}
 \newtheorem{example}[theorem]{Example}
 \newtheorem{question}[theorem]{Question}
 \newtheorem{hypothesis}{Hypothesis}
\newtheorem{case}{Case}
\newtheorem{exmp}[theorem]{Example}
  \newtheorem{rem}[theorem]{Remark}

\abovedisplayskip=3pt plus 1pt minus 1pt 
\belowdisplayskip=3pt plus 1pt minus 1pt 

\def\le{\leqslant}
\def\ge{\geqslant}
\def\dl{\displaystyle}

\numberwithin{equation}{section}




\vspace{8true mm}

\renewcommand{\baselinestretch}{1.9}\baselineskip 19pt

\noindent{\LARGE\bf BGP-Reflection Functors and Lusztig's Symmetries of Modified Quantized Enveloping Algebras}
\vspace{0.5 true cm}

\noindent{\normalsize\sf XIAO Jie$^{1}$ \& ZHAO
Minghui$^{1}$ \footnotetext{\baselineskip 10pt
This work was supported by NSF of China (No. 11131001).}}

\vspace{0.2 true cm}
\renewcommand{\baselinestretch}{1.5}\baselineskip 12pt
\noindent{\footnotesize\rm $^1$Department of Mathematical Sciences, Tsinghua University, Beijing 100084, P. R. China \\
(email: jxiao@math.tsinghua.edu.cn, zhaomh08@mails.tsinghua.edu.cn)\vspace{4mm}}

\baselineskip 12pt \renewcommand{\baselinestretch}{1.18}
\noindent{{\bf Abstract}\small\hspace{2.8mm} 
Let $\mathbf{U}$ be the quantized enveloping algebra and $\dot{\mathbf{U}}$ its modified form. Lusztig gives some symmetries on $\mathbf{U}$ and $\dot{\mathbf{U}}$.
Since the realization of $\mathbf{U}$ by the reduced Drinfeld double of the Ringel-Hall algebra, one can apply the BGP-reflection functors to the double Ringel-Hall algebra to obtain Lusztig's symmetries on $\mathbf{U}$ and their important properties, for instance, the braid relations.
In this paper, we define a modified form $\dot{\mathcal{H}}$ of the Ringel-Hall algebra and realize the Lusztig's symmetries on $\dot{\mathbf{U}}$
by applying the BGP-reflection functors to $\dot{\mathcal{H}}$.
 }

\noindent{\footnotesize{\bf Keywords:\hspace{2mm}
BGP-reflection functors, Lusztig's symmetries, Ringel-Hall algebras
}}


\section{Introduction}

Let $\mathbf{U}$ be the quantized enveloping algebra associated to a symmetrizable generalized Cartan matrix. Lusztig introduces some symmetries $T_i$ acting on an integrable $\mathbf{U}$-module and then on the quantized enveloping algebra $\mathbf{U}$ (\cite{Lusztig quantum deformation}\cite{Lusztig quantum groups at roots of 1}\cite{Lusztig introduction to quantum groups}). Let $\dot{\mathbf{U}}$ be the modified quantized enveloping algebra obtained from $\mathbf{U}$ by modifying the Cartan part $\mathbf{U}^{0}$ to $\oplus_{\lambda\in P}\mathbb{Q}(v)\mathbf{1}_{\lambda}$. This algebra has same representations with $\mathbf{U}$. Lusztig also introduces some symmetries $T_i$ acting on the modified quantized enveloping algebra $\dot{\mathbf{U}}$ (\cite{Lusztig introduction to quantum groups}).

Let $\mathcal{H}_q^{\ast}(\Lambda)$ be the Ringel-Hall algebra associated to a finite dimensional hereditary algebra $\Lambda.$  Then the composition subalgebra $\mathcal{C}_q^{\ast}(\Lambda)$ realizes the positive part $\mathbf{U}^+$ of the quantized enveloping algebra by the Ringel-Green Theorem (\cite{Green}\cite{Ringel Hall algebras and quantum groups}).  One can extend the Ringel-Green theorem to the Drinfeld double version and  realize the whole $\mathbf{U}$ by the reduced Drinfeld double of the composition algebra (\cite{X Drinfed Double and Ringel-Green Theorem and Hall Algebras}). These work give a connection between the representation theory of finite dimensional hereditary algebras and quantized enveloping algebras.

Via the Ringel-Hall algebra approach, one can apply the BGP-reflection functors to the quantum enveloping algebras $\mathbf{U}^+$ and $\mathbf{U}$ to obtain Lusztig's symmetries and their properties in a conceptual way (\cite{SV On the double of the Hall algebra of a quiver}\cite{XY}). This method gives a precise construction of Lusztig's symmetries not only in the quantum enveloping algebras, also for  the whole Drinfeld doubles of Ringel-Hall algebras (\cite{DX On double Ringel¨CHall algebras}\cite{DX Ringel¨CHall algebras and Lusztig's symmetries}).

In this paper, we define a modified form $\dot{\mathcal{H}}_q^{\ast}(\Lambda)$ of the Ringel-Hall algebra $\mathcal{H}_q^{\ast}(\Lambda)$.  We apply the BGP-reflection functors to obtain Lusztig's symmetries on  $\dot{\mathcal{H}}_q^{\ast}(\Lambda)$. Viewing the modified quantized enveloping algebra $\dot{\mathbf{U}}$ as a subalgebra of $\dot{\mathcal{H}}_q^{\ast}(\Lambda)$, we get a precise construction of Lusztig's symmetries on $\dot{\mathbf{U}}$. From this construction, we can obtain important properties of Lusztig's symmetries, for instance, the braid relations.

In Section 2, we first give the basic notation of quantized enveloping algebras and modified quantized enveloping algebras; then we recall the definition of Lusztig's symmetries on $\mathbf{U}$ and $\dot{\mathbf{U}}$. In Section 3, we recall the definition of the Ringel-Hall algebra $\mathcal{H}^{\ast}_q(\Lambda)$ and define a modified form $\dot{\mathcal{H}}_q^{\ast}(\Lambda)$ of it. In Section 4, we recall the BGP-reflection functors and define the corresponding maps from $\dot{\mathcal{H}}_q^{\ast}(\Lambda)$ to $\dot{\mathcal{H}}_q^{\ast}(\sigma_i\Lambda)$ induced by them. We prove in Section 6 that these maps induce algebra isomorphisms from $\dot{\mathbf{U}}$ to itself, which coincide to the Lusztig's symmetries on $\dot{\mathbf{U}}$ and satisfy the braid relations. In Section 5, we define Lusztig's symmetries on $\dot{\mathcal{H}}_q^{\ast}(\Lambda)$ and find the precise relation between these symmetries and the maps induced by the BGP-reflection functors.

\section{Quantized enveloping algebras and their modified forms}

\subsection{Quantized enveloping algebras}

Denote by $\mathbb{Q}$ the field of rational numbers and $\mathbb{Z}$ the ring of integers. Let $I$ be a finite index set with $|I|=n$ and $A=(a_{ij})_{i,j\in I}$ be a generalized Cartan matrix. Denote by $r(A)$ the rank of $A$. Let $P^{\vee}$ be a free abelian group of rank $2n-r(A)$ with a $\mathbb{Z}$-basis $\{h_i|i\in I\}\cup\{d_s|s=1,\ldots, n- r(A)\}$ and $\mathfrak{h}=\mathbb{Q}\otimes_{\mathbb{Z}}P^{\vee}$ be the $\mathbb{Q}$-linear space spanned by $P^{\vee}$. We call $P^{\vee}$ the dual weight lattice and $\mathfrak{h}$ the Cartan subalgebra. We also define the weight lattice to be $P=\{\lambda\in\mathfrak{h}^{\ast}|\lambda(P^{\vee})\subset\mathbb{Z}\}$.

Set $\Pi^{\vee}=\{h_i|i\in I\}$ and choose a linearly independent subset $\Pi=\{\alpha_i|i\in I\}\subset\mathfrak{h}^{\ast}$ satisfying $\alpha_j(h_i)=a_{ij}$ and $\alpha_j(d_s)=0$ or $1$ for $i,j\in I$, $s=1,\ldots, n-\textrm{rank} A$. The elements of $\Pi$ are called simple roots, and the elements of $\Pi^{\vee}$ are called simple coroots. The quintuple $(A,\Pi,\Pi^{\vee},P,P^{\vee})$ is called a Cartan datum associated with the generalized Cartan matrix $A$. Let $W$ be the Weyl group generated by simple reflections $s_i$ for all $i\in I$. There exists a bilinear form $(-,-)$ on $\mathfrak{h}^{\ast}$ (\cite{Kac}).

We recall the definition of the quantized enveloping algebras. Assume that  $A=(a_{ij})_{i,j\in I}$ is a symmetrizable generalized Cartan matrix and $D=\textrm{diag}(\varepsilon_i|i\in I)$ is its symmetrizing matrix.

Fix an indeterminate $v$. For $n\in \mathbb{Z}$, we set
\begin{displaymath}
[n]_{v}=\frac{v^n-v^{-n}}{v-v^{-1}},
\end{displaymath}
and $[0]_{v}!=1$,  $[n]_{v}!=[n]_{v}[n-1]_{v}\cdots[1]_{v}$ for $n\in\mathbb{Z}_{>0}$. For nonnegative integers $m\geq n\geq 0$, the analogues of binomial coefficients are given by
\begin{displaymath}
\left[{m \atop n}\right]_v=\frac{[m]_v!}{[n]_v![m-n]_v!}.
\end{displaymath}
Then $[n]_v$ and $\left[{m \atop n}\right]_v$ are elements of the field $\mathbb{Q}(v)$.

The quantized enveloping algebra $\mathbf{U}$ associated with a Cartan datum $(A,\Pi,\Pi^{\vee},P,P^{\vee})$ is an associative algebra over $\mathbb{Q}(v)$ with $\mathbf{1}$ generated by the elements $E_i$, $F_i(i\in I)$ and $K_{\mu}(\mu\in P^{\vee})$ subject to the following relations:
\begin{equation}
K_{0}=\mathbf{1}, K_{\mu}K_{\mu'}=K_{\mu+\mu'}\,\,\,\textrm{for all $\mu,\mu'\in P^{\vee}$};
\end{equation}
\begin{equation}
K_{\mu}E_{i}K_{-\mu}=v^{\alpha_i(\mu)}E_i\,\,\,\textrm{for all $i\in I$, $\mu\in P^{\vee}$};
\end{equation}
\begin{equation}
K_{\mu}F_{i}K_{-\mu}=v^{-\alpha_i(\mu)}E_i\,\,\,\textrm{for all $i\in I$, $\mu\in P^{\vee}$};
\end{equation}
\begin{equation}
E_iF_j-F_jE_i=\delta_{ij}\frac{{\tilde{K}}_{i}-{\tilde{K}}_{-i}}{v_i-v_i^{-1}}\,\,\,\textrm{for all $i,j\in I$};
\end{equation}
for $i\not=j$, setting $b=1-a_{ij}$,
\begin{equation}
\sum_{k=0}^{b}(-1)^{k}E_i^{(k)}E_jE_i^{(b-k)}=0;
\end{equation}
for $i\not=j$, setting $b=1-a_{ij}$,
\begin{equation}
\sum_{k=0}^{b}(-1)^{k}F_i^{(k)}F_jF_i^{(b-k)}=0.
\end{equation}
Here, $\tilde{K}_{\nu}=\Pi_{i\in I}K_{\varepsilon_i\nu_ih_i}$ for $\nu=\sum_{i\in I}\nu_ih_i$, $v_i=v^{\varepsilon_i}$ and $E_i^{(n)}=E_i^n/[n]_{v_i}!$, $F_i^{(n)}=F_i^n/[n]_{v_i}!$.

Let $\mathbf{U}^+$ (resp. $\mathbf{U}^-$) be the subalgebra of $\mathbf{U}$ generated by the elements $E_i$ (resp. $F_i$) for $i\in I$, and let $\mathbf{U}^{0}$ be the subalgebra of $\mathbf{U}$ generated by $K_{\mu}$ for $\mu\in P^{\vee}$. We know that the quantized enveloping algebra has the triangular decomposition
\begin{displaymath}
\mathbf{U}\cong {\mathbf{U}^-}\otimes{\mathbf{U}^{0}}\otimes{\mathbf{U}^{+}}.
\end{displaymath}

Let $\mathbf{f}$ be the associative algebra defined by Lusztig in \cite{Lusztig introduction to quantum groups}, which is generated by $\theta_i(i\in I)$ subject to the following relations
$$
\sum_{k=0}^{b}(-1)^{k}\theta_i^{(k)}\theta_j\theta_i^{(b-k)}=0,
$$
where $i\not=j$, $b=1-a_{ij}$ and $\theta_i^{(n)}=\theta_i^n/[n]_{v_i}!$.
There exist well-defined $\mathbb{Q}(v)$-algebra monomorphisms $\mathbf{f}\rightarrow \mathbf{U}(x\mapsto x^+)$ and $\mathbf{f}\rightarrow\mathbf{ U}(x\mapsto x^-)$ with image $\mathbf{U}^+$ and $\mathbf{U}^-$ respectively satisfying $E_i=\theta_i^+$ and $F_i=\theta_i^-$.



\subsection{Modified quantized enveloping algebras}

Let us recall the definition of the modified form $\dot{\mathbf{U}}$ of $\mathbf{U}$ in \cite{Lusztig introduction to quantum groups}.

If $\lambda',\lambda''\in P$, we set
\begin{displaymath}
_{\lambda'}\mathbf{U}_{\lambda''}=\mathbf{U}/\left(\sum_{\mu\in P^{\vee}}(K_{\mu}-v^{\lambda'(\mu)})\mathbf{U}+\sum_{\mu\in P^{\vee}}\mathbf{U}(K_{\mu}-v^{\lambda''(\mu)})\right).
\end{displaymath}
Let $\pi_{\lambda',\lambda''}:\mathbf{U}\rightarrow _{\lambda'}\mathbf{U}_{\lambda''}$ be the canonical projection and
\begin{displaymath}
\dot{\mathbf{U}}=\bigoplus_{\lambda',\lambda''\in P} {_{\lambda'}\mathbf{U}_{\lambda''}}.
\end{displaymath}

Consider the weight space decomposition $\mathbf{U}=\oplus_{\beta}\mathbf{U}(\beta)$, where $\beta$ runs through $\mathbb{Z}I$ and
$\mathbf{U}(\beta)=\{x\in\mathbf{U}|K_{\mu}xK_{\mu}^{-1}=v^{\beta(\mu)}x\,\,\,\textrm{for all $\mu\in P^{\vee}$}\}$.
The image of summands $\mathbf{U}(\beta)$ under $\pi_{\lambda',\lambda''}$ form the weight space decomposition $_{\lambda'}\mathbf{U}_{\lambda''}=\oplus_{\beta}{_{\lambda'}\mathbf{U}_{\lambda''}}(\beta)$.
Note that $_{\lambda'}\mathbf{U}_{\lambda''}(\beta)=0$ unless $\lambda'-\lambda''=\beta$.

There is a natural associative $\mathbb{Q}(v)$-algebra structure on $\dot{\mathbf{U}}$ inherited from that of $\mathbf{U}$. It is defined as follows: for any $\lambda'_1,\lambda''_1,\lambda'_2,\lambda''_2\in P$, $\beta_1,\beta_2\in \mathbb{Z}I$ such that $\lambda'_1-\lambda''_1=\beta_1,\lambda'_2-\lambda''_2=\beta_2$ and any $x\in \mathbf{U}(\beta_1),y\in \mathbf{U}(\beta_2)$,
\begin{displaymath}
\pi_{\lambda'_1,\lambda''_1}(x)\pi_{\lambda'_2,\lambda''_2}(y)=
\left\{
  \begin{array}{cc}
    \pi_{\lambda'_1,\lambda''_2}(xy) & \textrm{if}\,\,\,\lambda''_1=\lambda'_2\\
     0 & \textrm{otherwise}
  \end{array}
\right..
\end{displaymath}

Let $\mathbf{1}_{\lambda}=\pi_{\lambda,\lambda}(\mathbf{1})$, where $\mathbf{1}$ is the unit element of $\mathbf{U}$. Then they satisfy $\mathbf{1}_{\lambda}\mathbf{1}_{\lambda'}=\delta_{\lambda,\lambda'}\mathbf{1}_{\lambda}$. In general, there is no unit element in the algebra $\dot{\mathbf{U}}$. However the family $(\mathbf{1}_{\lambda})_{\lambda\in P}$ can be regarded locally as the unit element in $\dot{\mathbf{U}}$.

Note that ${_{\lambda'}\mathbf{U}_{\lambda''}}=\mathbf{1}_{\lambda'}\dot{\mathbf{U}}\mathbf{1}_{\lambda''}$. We define $\dot{\mathbf{U}}\mathbf{1}_{\lambda}=\oplus_{\lambda'\in{P}}\mathbf{1}_{\lambda'}\dot{\mathbf{U}}\mathbf{1}_{\lambda}$. Then $\dot{\mathbf{U}}=\oplus_{\lambda\in{P}}\dot{\mathbf{U}}\mathbf{1}_{\lambda}$.

\subsection{Lusztig's symmetries on $\dot{\mathbf{U}}$}

In \cite{Lusztig introduction to quantum groups}, Lusztig introduces some symmetries on $\mathbf{U}$, which is now called Lusztig's symmetries.

Fix $i\in I$. Define $T_i:\mathbf{U}\rightarrow \mathbf{U}$ on the generators as follows:
\begin{eqnarray*}
&&T_i(E_i)=-F_i\tilde{K}_{i},T_i(F_i)=-\tilde{K}_{-i}E_i;\\
&&T_i(E_j)=\sum_{r+s=-\alpha_j(h_i)}(-1)^rv_i^{-r}E_i^{(s)}E_jE_i^{(r)} \textrm{for $j\neq i$};\\
&&T_i(F_j)=\sum_{r+s=-\alpha_j(h_i)}(-1)^rv_i^{r}F_i^{(r)}F_jF_i^{(s)} \textrm{for $j\neq i$};\\
&&T_i(K_{\mu})=K_{\mu-\alpha_{i}(\mu)h_i}.
\end{eqnarray*}

Lusztig also introduces symmetries $T_i:\dot{\mathbf{U}}\rightarrow \dot{\mathbf{U}}$ induced by the symmetries on $\mathbf{U}$. We write the following formulas:
\begin{eqnarray*}
&&T_i(E_i\mathbf{1}_{\lambda})=-v_i^{-\lambda(h_i)}F_i\mathbf{1}_{s_i\lambda};\\
&&T_i(F_i\mathbf{1}_{\lambda})=-v_i^{-(2-\lambda(h_i))}E_i\mathbf{1}_{s_i\lambda};\\
&&T_i(E_j\mathbf{1}_{\lambda})=\sum_{r+s=-\alpha_j(h_i)}(-1)^rv_i^{-r}E_i^{(s)}E_jE_i^{(r)}\mathbf{1}_{s_i\lambda}\,\,\,\textrm{for $j\neq i$};\\
&&T_i(F_j\mathbf{1}_{\lambda})=\sum_{r+s=-\alpha_j(h_i)}(-1)^rv_i^{r}F_i^{(r)}F_jF_i^{(s)}\mathbf{1}_{s_i\lambda}\,\,\,\textrm{for $j\neq i$}.
\end{eqnarray*}

\section{Ringel-Hall algebras and their modified form}

\subsection{Ringel-Hall algebras}

In this subsection, we recall the definition of Ringel-Hall algebras, following the notations in \cite{Dlab Ringel}, \cite{XY} and \cite{DX Ringel¨CHall algebras and Lusztig's symmetries}.

Let $k$ be a finite field and $\Lambda$ be a finite dimensional hereditary $k$-algebra. According to \cite{Dlab Ringel}, we can identity $\Lambda$ with the tensor algebra of a $k$-species. A valued graph $(\Gamma,\mathbf{d})$ is a finite set $\Gamma$ together with nonnegative integers $d_{ij}$ for all $i,j\in\Gamma$ such that $d_{ii}=0$ and there exist positive integers $\{\varepsilon_i\}_{i\in\Gamma}$ satisfying $$d_{ij}\varepsilon_j=d_{ji}\varepsilon_i\,\,\,\,\textrm{for $i,j\in\Gamma$}.$$
Given a Cartan datum $(A,\Pi,\Pi^{\vee},P,P^{\vee})$, there is a valued graph $(\Gamma,\mathbf{d})$ corresponding to it.

An orientation $\Omega$ of a valued graph $(\Gamma,\mathbf{d})$ is given by an order on each edge $\{i,j\}$, which is indicated by an arrow $i\rightarrow j$. We call $Q=(\Gamma,\mathbf{d},\Omega)$ a valued quiver.

We assume that $Q=(\Gamma,\mathbf{d},\Omega)$ is connected and contains no cycles. Let $\mathcal{S}=(F_i,_iM_j)_{i,j\in\Gamma}$ be a reduced $k$-species of type $Q$, that is, for all $i,j\in\Gamma$, $_iM_j$ is an $F_i$-$F_j$-bimodule, where $F_i$ and $F_j$ are finite extensions of $k$ in an algebraic closure and $\dim(_iM_j)_{F_j}=d_{ij}$ and $\dim_k(F_i)=\varepsilon_i$. A $k$-representation $(V_i,_j\varphi_i)$ of $\mathcal{S}$ is given by vector spaces $(V_i)_{F_i}$ for any $i\in\Gamma$ and $F_j$-linear mapping $_j\varphi_i:V_i\otimes_iM_j\rightarrow V_j$ for any $i\rightarrow j$. Such a representation is called finite dimensional if $\sum_{i\in\Gamma}\dim_kV_i<\infty$. We denote by rep-$\mathcal{S}$ the category of finite dimensional representations of $\mathcal{S}$ over $k$. Let $\Lambda$ be the tensor algebra of $\mathcal{S}$. Then the category rep-$\mathcal{S}$ is equivalent to the module category mod-$\Lambda$ of finite dimensional modules over $\Lambda$.

Given three modules $L,M$ and $N$ in mod-$\Lambda$, denote by $g^{L}_{MN}$ the number of $\Lambda$-submodules $W$ of $L$ such that $W\simeq N$ and $L/W\simeq N$ in mod-$\Lambda$. Let $v=\sqrt{|k|}\in \mathbb{C}$, $\mathcal{P}$ be the set of isomorphism classes of finite dimensional (nilpotent) $\Lambda$-modules and ind$(\mathcal{P})$ be the set of isomorphism classes of indecomposable finite dimensional (nilpotent) $\Lambda$-modules. The Ringel-Hall algebra $\mathcal{H}_q(\Lambda)$ of $\Lambda$ is by definition the $\mathbb{Q}(v)$-space with basis $\{u_{[M]}|[M]\in\mathcal{P}\}$ whose multiplication is given by
\begin{displaymath}
u_{[M]}u_{[N]}=\sum_{[L]\in\mathcal{P}}g^{L}_{MN}u_{[L]}.
\end{displaymath}
It is easily seen that $\mathcal{H}_q(\Lambda)$ is associative $\mathbb{Q}(v)$-algebra with unit $u_{[0]}$, where $0$ denotes the zero module.

For each representation $V=(V_i,_j\varphi_i)$ in rep-$\mathcal{S}$, the dimension vector of $V$ is defined to be $\underline{\dim}V=(\dim_{F_i}V_i)_{i\in\Gamma}\in\mathbb{N}^{\Gamma}$. For $V,W\in\textrm{rep-$\mathcal{S}$}$, The Euler form is defined by
$$\langle\underline{\dim}V,\underline{\dim}W\rangle=\sum_{i\in\Gamma}\varepsilon_ia_ib_i-\sum_{i\rightarrow j}d_{ij}\varepsilon_ja_ib_j,$$
where $\underline{\dim}V=(a_1,\ldots,a_n)$ and $\underline{\dim}W=(b_1,\ldots,b_n)$. It is well known that
$$\langle\underline{\dim}V,\underline{\dim}W\rangle=\dim_{k}\textrm{Hom}_{\Lambda}(V,W)-\dim_{k}\textrm{Ext}_{\Lambda}(V,W).$$
Further, the symmetric Euler form is defined as
$$(\underline{\dim}V,\underline{\dim}W)=\langle\underline{\dim}V,\underline{\dim}W\rangle+\langle\underline{\dim}W,\underline{\dim}V\rangle.$$
Both $\langle-,-\rangle$ and $(-,-)$ are well defined on the Grothendieck group $G(\Lambda)$ of mod-$\Lambda$. In fact, the Grothendieck group $G(\Lambda)$ with the symmetric Euler form is a Cartan datum.

Let $I\subset\mathcal{P}$ be the set of isomorphism classes of (nilpotent) simple $\Lambda$-modules, which can be identified with $\Gamma$. Then the Euler form and the symmetric Euler form are defined on $\mathbb{Z}I$. We also identify $\mathbb{N}^{\Gamma}$ with $\mathbb{N}I$ and regard $\underline{\dim}V$ as an element in $\mathbb{N}I$ for each representation $V=(V_i,_j\varphi_i)$ in rep-$\mathcal{S}$. For each $\alpha\in\mathcal{P}$, we fix a representation $V_{\alpha}$ in the isomorphism class $\alpha$. For $\alpha,\beta\in\mathcal{P}$, we set
$$\langle\alpha,\beta\rangle=\langle\underline{\dim}V_{\alpha},\underline{\dim}V_{\beta}\rangle$$
and
$$(\alpha,\beta)=(\underline{\dim}V_{\alpha},\underline{\dim}V_{\beta}).$$
Note that for $\alpha,\beta\in\mathcal{P}$, $(\alpha,\beta)=(\sum_{i\in I}a_i\alpha_i,\sum_{i\in I}b_i\alpha_i)$, where $\underline{\dim} V_{\alpha}=\sum a_ii$ and $\underline{\dim} V_{\beta}=\sum b_ii$. Hence, we also use $\alpha$ to express the element $\sum_{i\in I}a_i\alpha_i$ in $P$ and the element $\sum_{i\in I}a_ih_i$ in $P^{\vee}$.

The twisted Ringel-Hall algebra $\mathcal{H}^{\ast}_q(\Lambda)$ is defined as follows. Set $\mathcal{H}^{\ast}_q(\Lambda)=\mathcal{H}_q(\Lambda)$ as $\mathbb{Q}(v)$-vector space and define the multiplication by
\begin{displaymath}
u_{[M]}\ast u_{[N]}=v^{\langle\underline{\dim}M,\underline{\dim}N\rangle}\sum_{[L]\in\mathcal{P}}g^{L}_{MN}u_{[L]}.
\end{displaymath}
The composition algebra $\mathcal{C}^{\ast}_q(\Lambda)$ is a subalgebra of $\mathcal{H}^{\ast}_q(\Lambda)$ generated by $u_i=u_{[S_i]}$, $i\in I$, where $S_i$ is the (nilpotent) simple module corresponding to $i\in I$.
For any $\Lambda$-module $M$, we denote $$\langle M\rangle=v^{-\dim M+\dim \textrm{End}_{\Lambda}(M)}u_{[M]}.$$
Note that $\{\langle M\rangle|M\in\mathcal{P}\}$ is a $\mathbb{Q}(v)$-basis of $\mathcal{H}^{\ast}_q(\Lambda)$.



Then we consider the generic form of Ringel-Hall algebras. Let $Q$ be a valued quiver and $\Lambda_k$ the corresponding finite dimensional hereditary algebra of a $k$-species which is of type $Q$. Denote by $\mathcal{H}_q^{\ast}(\Lambda_k)$ the twisted Ringel-Hall algebra of $\Lambda_k$. Let $\mathcal{K}$ be a set of finite fields $k$ such that the set $\{q_k=|k||k\in\mathcal{K}\}$ is infinite and ${R}$ be an integral domain containing $\mathbb{Q}$ and an element $v_{q_k}$ such that $v_{q_k}^2=q_k$ for each $k\in\mathcal{K}$. For each $k\in\mathcal{K}$, we consider the composition algebra $\mathcal{C}^{\ast}_q(\Lambda_k)$ which is the ${R}$-subalgebra of $\mathcal{H}_q^{\ast}(\Lambda_k)$ generated by the elements $u_i(k)$. Consider the direct product
\begin{displaymath}
\mathcal{H}^{\ast}(Q)=\prod_{k\in\mathcal{K}}\mathcal{H}_q^{\ast}(\Lambda_k)
\end{displaymath}
and the elements $v=(v_{q_k})_{k\in\mathcal{K}}$, $v^{-1}=(v_{q_k}^{-1})_{k\in\mathcal{K}}$ and $u_i=(u_i(k))_{k\in\mathcal{K}}$. By $\mathcal{C}^{\ast}(Q)_{\mathcal{A}}$ we denote the subalgebra of $\mathcal{H}^{\ast}(Q)$ generated by $v$, $v^{-1}$ and $u_i$ over  $\mathbb{Q}$, where $\mathcal{A}=\mathbb{Q}[v,v^{-1}]$. We may regard it as the $\mathcal{A}$-algebra generated by $u_i$ where $v$ is considered as an indeterminate. Finally, denote by $\mathcal{C}^{\ast}(Q)=\mathbb{Q}(v)\otimes\mathcal{C}^{\ast}(Q)_{\mathcal{A}}$ the generic twisted composition algebra of type $Q$.

\begin{remark}
If $Q$ is a Dynkin quiver, then the generic composition algebra of $Q$ can be defined directly using Hall polynomials.
\end{remark}

Then we have the following well-known result of Green and Ringel (\cite{Green}\cite{Ringel Hall algebras and quantum groups}).

\begin{theorem}\label{theorem1}
Let $Q$ be a valued quiver, $A$ be the associated  generalized Cartan matrix, and $\mathbf{f}$ be the Lusztig's algebra of type $A$. Then the correspondence $u_i\mapsto\theta_i$, $i\in I$ induces an algebra isomorphism from $\mathcal{C}^{\ast}(Q)$ to $\mathbf{f}$.
\end{theorem}

\subsection{Double Ringel-Hall algebras}

Let $\Lambda$ be a finite dimensional hereditary algebra. In \cite{X Drinfed Double and Ringel-Green Theorem and Hall Algebras}, the reduced Drinfeld double $\mathcal{D}(\Lambda)$ of $\Lambda$ is defined. As an associative algebra, $\mathcal{D}(\Lambda)$ is generated by $\langle{u_{\alpha}(+)}\rangle$, $\langle{u_{\alpha}(-)}\rangle(\alpha\in\mathcal{P})$ and $K_{\mu}(\mu\in P^{\vee})$ subject to the following relations (\cite{XY}):
\begin{eqnarray}
&&K_0=\langle{u_0(+)}\rangle=\langle{u_0(-)}\rangle=\mathbf{1},\,\,K_{\mu}K_{\mu'}=K_{\mu+\mu'};\\
&&\langle{u_{\alpha}(+)}\rangle\langle{u_{\beta}(+)}\rangle
=v^{-\langle{\beta,\alpha}\rangle}\sum_{\lambda\in\mathcal{P}}g_{\alpha\beta}^{\lambda}\langle{u_{\lambda}(+)}\rangle;\\
&&\langle{u_{\alpha}(-)}\rangle\langle{u_{\beta}(-)}\rangle
=v^{-\langle{\beta,\alpha}\rangle}\sum_{\lambda\in\mathcal{P}}g_{\alpha\beta}^{\lambda}\langle{u_{\lambda}(-)}\rangle;\\
&&K_{\mu}\langle{u_{\beta}(+)}\rangle=v^{\beta(\mu)}\langle{u_{\beta}(+)}\rangle K_{\mu};\\
&&K_{\mu}\langle{u_{\beta}(-)}\rangle=v^{-\beta(\mu)}\langle{u_{\beta}(-)}\rangle K_{\mu};\\
&&\sum_{\alpha,\alpha'\in\mathcal{P}}v^{\langle{\alpha',\alpha}\rangle+(\alpha,\alpha)}\frac{a_{\alpha'}}{a_{\lambda'}}g_{\alpha'\alpha}^{\lambda'}\tilde{K}_{-\alpha}\langle{u_{\alpha'}(-)}\rangle r'_{\alpha}(\langle{u_{\lambda}(+)}\rangle)\nonumber\\
&=&\sum_{\alpha,\beta\in\mathcal{P}}v^{\langle{\alpha,\beta}\rangle+(\beta,\beta)}\frac{a_{\alpha}}{a_{\lambda}}g_{\alpha\beta}^{\lambda}\tilde{K}_{\beta}\langle{u_{\alpha}(+)}\rangle r_{\beta}(\langle{u_{\lambda'}(-)}\rangle),
\end{eqnarray}
where $\alpha,\beta,\lambda,\lambda'\in\mathcal{P}$, $\mu,\mu'\in P^{\vee}$ and
\begin{displaymath}
r'_{\alpha}(\langle {u_{\lambda}(+)}\rangle)=\sum_{\beta\in\mathcal{P}}v^{\langle{\alpha,\beta}\rangle+(\alpha,\beta)}g_{\alpha\beta}^{\lambda}\frac{a_{\alpha}a_{\beta}}{a_{\lambda}}\langle u_{\beta}(+)\rangle;
\end{displaymath}
\begin{displaymath}
r_{\alpha}(\langle {u_{\lambda}(-)}\rangle)=\sum_{\beta\in\mathcal{P}}v^{\langle{\alpha,\beta}\rangle+(\alpha,\beta)}g_{\alpha\beta}^{\lambda}\frac{a_{\alpha}a_{\beta}}{a_{\lambda}}\langle u_{\beta}(-)\rangle.
\end{displaymath}

From the definition of $\mathcal{D}(\Lambda)$, we have two algebra monomorphisms $(+):\mathcal{H}_q^{\ast}(\Lambda)\rightarrow\mathcal{D}(\Lambda)$ mapping $\langle{M(\lambda)}\rangle$ to $u_{\lambda}(+)$ and $(-):\mathcal{H}_q^{\ast}(\Lambda)\rightarrow\mathcal{D}(\Lambda)$ mapping $\langle{M(\lambda)}\rangle$ to $u_{\lambda}(-)$ for all $\lambda\in\mathcal{P}$.

Consider the weight space decomposition $\mathcal{D}(\Lambda)=\oplus_{\beta}\mathcal{D}(\Lambda)(\beta)$, where $\beta$ runs through $\mathbb{Z}I$ and
$\mathcal{D}(\Lambda)(\beta)=\{x\in\mathcal{D}(\Lambda)|K_{\mu}xK_{\mu}^{-1}=v^{\beta(\mu)}x\,\,\,\textrm{for all $\mu\in P^{\vee}$}\}$.

Let $\mathcal{D}_c(\Lambda)$ be the subalgebra of $\mathcal{D}(\Lambda)$ generated by $\langle{u_i(\pm)}\rangle(i\in I)$ and $K_{\mu}(\mu\in P^{\vee})$. In \cite{X Drinfed Double and Ringel-Green Theorem and Hall Algebras}, the Green-Ringel Theorem \ref{theorem1} is extended to the Drinfeld double version and $\mathcal{D}_c(\Lambda)$ realizes the corresponding quantum enveloping algebra $\mathbf{U}$.

\subsection{Another definition of $\dot{\mathbf{U}}$ and a similar form of $\mathcal{H}^{\ast}(\Lambda)$}

In \cite{Lusztig introduction to quantum groups}, Lusztig gives another definition of $\dot{\mathbf{U}}$ as follows. $\dot{\mathbf{U}}$ can be viewed as the algebra generated by the symbols $x^+\mathbf{1}_{\zeta}x'^-$ and $x^-\mathbf{1}_{\zeta}x'^+$ with $x\in\mathbf{f}_{\nu},x'\in\mathbf{f}_{\nu'}$ for various $\nu,\nu'\in\mathbb{N}I$ and $\zeta\in P$; these symbols are subject to the following relations (\ref{3}) to (\ref{4}):
\begin{equation}
(\theta_i^{(a)})^+\mathbf{1}_{\zeta}(\theta_j^{(b)})^-=(\theta_j^{(b)})^-\mathbf{1}_{\zeta+a\alpha_i+b\alpha_j}(\theta_i^{(a)})^+ \textrm{if $i\neq j$};\label{3}
\end{equation}
\begin{equation}
(\theta_i^{(a)})^+\mathbf{1}_{-\zeta}(\theta_i^{(b)})^-=\sum_{t\geq 0}\left[\begin{array}{c}
                                                                     a+b-\zeta(h_i) \\
                                                                     t
                                                                   \end{array}
\right]_{v_i}(\theta_i^{(b-t)})^-\mathbf{1}_{-\zeta+(a+b-t)\alpha_i}(\theta_i^{(a-t)})^+;
\end{equation}
\begin{equation}
(\theta_i^{(b)})^-\mathbf{1}_{\zeta}(\theta_i^{(a)})^+=\sum_{t\geq 0}\left[\begin{array}{c}
                                                                     a+b-\zeta(h_i) \\
                                                                     t
                                                                   \end{array}
\right]_{v_i}(\theta_i^{(a-t)})^+\mathbf{1}_{\zeta-(a+b-t)\alpha_i}(\theta_i^{(b-t)})^+;
\end{equation}
\begin{equation}
x^+\mathbf{1}_{\zeta}=\mathbf{1}_{\zeta+\nu}x^+,x^-\mathbf{1}_{\zeta}=\mathbf{1}_{\zeta-\nu}x^- \textrm{for $x\in\mathbf{f}_{\nu}$};
\end{equation}
\begin{equation}
(x^+\mathbf{1}_{\zeta})(\mathbf{1}_{\zeta'}x'^-)=\delta_{\zeta,\zeta'}x^+\mathbf{1}_{\zeta}x'^-,
(x^-\mathbf{1}_{\zeta})(\mathbf{1}_{\zeta'}x'^+)=\delta_{\zeta,\zeta'}x^-\mathbf{1}_{\zeta}x'^+;
\end{equation}
\begin{eqnarray}
&&(x^+\mathbf{1}_{\zeta})(\mathbf{1}_{\zeta'}x'^+)=\delta_{\zeta,\zeta'}\mathbf{1}_{\zeta+\nu}(xx')^+,\nonumber\\
&&(x^-\mathbf{1}_{\zeta})(\mathbf{1}_{\zeta'}x'^-)=\delta_{\zeta,\zeta'}\mathbf{1}_{\zeta-\nu}(xx')^-\textrm{for $x\in\mathbf{f}_{\nu}$};
\end{eqnarray}
\begin{eqnarray}
&&(rx+r'x')^+\mathbf{1}_{\zeta}=rx^+\mathbf{1}_{\zeta}+r'x'^+\mathbf{1}_{\zeta},(rx+r'x')^-\mathbf{1}_{\zeta}=rx^-\mathbf{1}_{\zeta}+r'x'^-\mathbf{1}_{\zeta}\nonumber\\
&&\textrm{for $x,x'\in\mathbf{f}_{\nu}$ and $r,r'\in\mathbb{Q}(v)$}.\label{4}
\end{eqnarray}

Let $k$ be a finite field and $\Lambda$ a finite dimensional hereditary $k$-algebra. For each $\nu\in\mathbb{N}I$, set $$\mathcal{H}_q^{\ast}(\Lambda)_{\nu}=\textrm{span}\{u_{[M]}|\underline{\dim} M=\nu\}.$$
Similarly, we can define $\dot{\mathcal{H}}_q^{\ast}(\Lambda)$ as follows.
$\dot{\mathcal{H}}_q^{\ast}(\Lambda)$ is the algebra generated by the symbols $x^+\mathbf{1}_{\zeta}x'^-$ and $x^-\mathbf{1}_{\zeta}x'^+$
with $x\in\mathcal{H}_q^{\ast}(\Lambda)_{\nu},x'\in\mathcal{H}_q^{\ast}(\Lambda)_{\nu'}$ for various $\nu,\nu'\in\mathbb{N}I$ and $\zeta\in P$; these symbols are subject to the following relations (\ref{1}) to (\ref{2}):

\begin{eqnarray}
&&\sum_{\alpha,\alpha'\in\mathcal{P}}v^{\langle\alpha',\alpha\rangle+(\alpha,\alpha)+(\zeta,-\alpha)}\frac{a_{\alpha'}}{a_{\lambda'}}g_{\alpha'\alpha}^{\lambda'}(-1)^{tr\alpha'}v^{m(\alpha')}\langle M(\alpha')\rangle^-\mathbf{1}_{\zeta+\alpha'}(r'_{\alpha}(\langle M(\lambda)\rangle))^+\nonumber\\
&&=\sum_{\alpha,\beta\in\mathcal{P}}v^{\langle\alpha,\beta\rangle+(\beta,\beta)+(\zeta,\beta)}\frac{a_{\alpha}}{a_{\lambda}}g_{\alpha\beta}^{\lambda}(-1)^{tr(\lambda'-\beta)}v^{m(\lambda'-\beta)}\langle M(\alpha)\rangle^+\mathbf{1}_{\zeta-\alpha}(r_{\beta}(\langle M(\lambda')\rangle))^-;\label{1}
\end{eqnarray}
\begin{equation}
x^+\mathbf{1}_{\zeta}=\mathbf{1}_{\zeta+\nu}x^+,x^-\mathbf{1}_{\zeta}=\mathbf{1}_{\zeta-\nu}x^- \textrm{for $x\in\mathcal{H}_q^{\ast}(\Lambda)_{\nu}$};
\end{equation}
\begin{equation}
(x^+\mathbf{1}_{\zeta})(\mathbf{1}_{\zeta'}x'^-)=\delta_{\zeta,\zeta'}x^+\mathbf{1}_{\zeta}x'^-,(x^-\mathbf{1}_{\zeta})(\mathbf{1}_{\zeta'}x'^+)=\delta_{\zeta,\zeta'}x^-\mathbf{1}_{\zeta}x'^+;
\end{equation}
\begin{eqnarray}
&&(x^+\mathbf{1}_{\zeta})(\mathbf{1}_{\zeta'}x'^+)=\delta_{\zeta,\zeta'}\mathbf{1}_{\zeta+\nu}(xx')^+,\nonumber\\
&&(x^-\mathbf{1}_{\zeta})(\mathbf{1}_{\zeta'}x'^-)=\delta_{\zeta,\zeta'}\mathbf{1}_{\zeta-\nu}(xx')^-\textrm{for $x\in\mathcal{H}_q^{\ast}(\Lambda)_{\nu}$};\label{5}
\end{eqnarray}
\begin{eqnarray}
&&(rx+r'x')^+\mathbf{1}_{\zeta}=rx^+\mathbf{1}_{\zeta}+r'x'^+\mathbf{1}_{\zeta},(rx+r'x')^-\mathbf{1}_{\zeta}=rx^-\mathbf{1}_{\zeta}+r'x'^-\mathbf{1}_{\zeta}\nonumber\\
&&\textrm{for $x,x'\in\mathcal{H}_q^{\ast}(Q)_{\nu}$ and $r,r'\in\mathbb{Q}(v)$}.\label{2}
\end{eqnarray}
Here $a_{\lambda}$ is the order of the automorphism group of $V_{\lambda}$ for $\lambda\in\mathcal{P}$, $tr\alpha=\sum_{i\in I}a_i$, $m(\alpha)=\sum_{i\in I}a_i\varepsilon_i$ if $\alpha=\sum_{i\in I}a_i\alpha_i$,
and
\begin{displaymath}
r_{\alpha}(\langle M(\lambda)\rangle)=\sum_{\beta\in\mathcal{P}}v^{\langle{\beta,\alpha}\rangle+(\beta,\alpha)}g_{\beta\alpha}^{\lambda}\frac{a_{\beta}a_{\alpha}}{a_{\lambda}}\langle M(\beta)\rangle;
\end{displaymath}
\begin{displaymath}
r'_{\alpha}(\langle M(\lambda)\rangle)=\sum_{\beta\in\mathcal{P}}v^{\langle{\alpha,\beta}\rangle+(\beta,\alpha)}g_{\alpha\beta}^{\lambda}\frac{a_{\alpha}a_{\beta}}{a_{\lambda}}\langle M(\beta)\rangle.
\end{displaymath}

Similarly to the case of modified form of quantum group, we have the following direct sums decompositions
$$\dot{\mathcal{H}}_q^{\ast}(\Lambda)=\bigoplus_{\zeta\in P}\{x^+1_{\zeta}x'^-|x,x'\in\mathcal{H}_q^{\ast}(\Lambda)\}$$
and
$$\dot{\mathcal{H}}_q^{\ast}(\Lambda)=\bigoplus_{\zeta\in P}\{x^-1_{\zeta}x'^+|x,x'\in\mathcal{H}_q^{\ast}(\Lambda)\}.$$

Let $\dot{\mathcal{C}}_q^{\ast}(\Lambda)$ be the composition algebra, which is a subalgebra of $\dot{\mathcal{H}}_q^{\ast}(\Lambda)$ generated by $u_i^+1_{\zeta}u_j^-$ and $u_i^-1_{\zeta}u_j^+$ for all $i,j\in I$ and $\zeta\in P$.

Similarly to the Ringel-Hall algebra case we can consider the generic form of $\dot{\mathcal{H}}^{\ast}(\Lambda)$ and its generic composition subalgebra $\dot{\mathcal{C}}^{\ast}(Q)$, which is isomorphic to the corresponding modified quantum enveloping algebra $\dot{\mathbf{U}}$. If a formula in $\dot{\mathcal{C}}_q^{\ast}(\Lambda)$ is independent of the choice of the field, it can be viewed as a formula in $\dot{\mathcal{C}}^{\ast}(Q)\simeq\dot{\mathbf{U}}$.

\section{BGP-reflection functors and Lusztig's symmetries}

In this section we apply the BGP-reflection functors to the Ringel-Hall algebras and obtain an alternative construction of Lusztig's symmetries on modified quantum enveloping algebras.

\subsection{BGP-reflection functors}

Let $Q=(\Gamma,\mathbf{d},\Omega)$ be a valued quiver, $\mathcal{S}=(F_i,_iM_j)_{i,j\in\Gamma}$ be a $k$-species of type $Q$ and $p$ be a sink or source of $(\Gamma,\mathbf{d},\Omega)$. We define a new orientation $\sigma_p\Omega$ of $(\Gamma,\mathbf{d})$ by reversing the direction of arrows along all edges containing $p$ and $\sigma_pQ=(\Gamma,\mathbf{d},\sigma_p\Omega)$. Let $\sigma_p\mathcal{S}$ be the $k$-species obtained from $\mathcal{S}$ by replacing $_rM_s$ by its $k$-dual for $r=p$ or $s=p$. Then $\sigma_p\mathcal{S}$ is a reduced $k$-species of type $\sigma_pQ$. Assume $\Lambda$ is the corresponding finite dimensional hereditary algebra to $\mathcal{S}$. We denote by $\sigma_i\Lambda$ the corresponding finite dimensional hereditary algebra to $\sigma_i\mathcal{S}$.

Now, we recall the definition of the Bernstein-Gelfand-Ponomarev (BGP) reflection functors
$\sigma_p^{\pm}:\textrm{rep-$\mathcal{S}\rightarrow$rep-$\sigma_p\mathcal{S}$}$ (\cite{BGP} \cite{Dlab Ringel} \cite{XY}).

Let $p$ be a sink of $\Omega$. For any $V=(V_i,_j\varphi_i)\in\textrm{rep-}\mathcal{S}$, define $\sigma_p^{+}V=W=(W_i,_j\psi_i)$ as follows.
Let $$W_i=V_i\,\,\,\,\,\textrm{for $i\neq p$},$$
and $W_p$ be the kernel of
$$\xymatrix{\bigoplus_{j\rightarrow p}V_j\otimes_jM_p \ar[r]^-{(_p\varphi_j)_j}&V_p},$$
that is, we have the following exact sequence of vector spaces
$$\xymatrix{0\ar[r]&W_p\ar[r]^-{(_j\kappa_p)_j}&\bigoplus_{j\rightarrow p}V_j\otimes_jM_p \ar[r]^-{(_p\varphi_j)_j}&V_p}.$$
Let $$_j\psi_i=_j\varphi_i\,\,\,\,\,\textrm{for $i\neq p$},$$
and
$$_j\psi_p=_j\bar{\kappa}_p:W_p\otimes_pM_j\rightarrow W_j,$$
where $_j\bar{\kappa}_p$ corresponds to $_j\kappa_p$ under the natural isomorphism
$$\textrm{Hom}_{F_j}(W_p\otimes_pM_j,W_j)\simeq \textrm{Hom}_{F_p}(W_p,W_j\otimes_jM_p).$$
For any morphism $f=(f_i):V\rightarrow V'$ in rep-$\mathcal{S}$, define $\sigma_p^+f=g=(g_i)$ as follows.
Let $$g_i=h_i\,\,\,\,\,\textrm{for $i\neq p$}$$
and $g_p:W_p\rightarrow W'_p$ be the restriction of $\oplus_{j\rightarrow p}(f_j\otimes 1)$, that is, we have the following commutative diagram
$$\xymatrix{0\ar[r]&W_p\ar[d]^-{g_p}\ar[r]^-{(_j\kappa_p)_j}&\bigoplus_{j\rightarrow p}V_j\otimes_jM_p \ar[d]^-{\bigoplus_{j\rightarrow p}(f_j\otimes 1)}\ar[r]^-{(_p\varphi_j)_j}&V_p\ar[d]^-{f_p}\\
0\ar[r]&W'_p\ar[r]^-{(_j\kappa'_p)_j}&\bigoplus_{j\rightarrow p}V'_j\otimes_jM_p \ar[r]^-{(_p\varphi'_j)_j}&V'_p}$$

%
Similarly, if $p$ is a source of $\Omega$, we can define $\sigma_p^-$ from rep-$\mathcal{S}$ to rep-$\sigma_p\mathcal{S}$.

For $i\in\Gamma$, let rep-$\mathcal{S}\langle i\rangle$ be the full subcategory of rep-$\mathcal{S}$ containing all representations which do not have $V_i$ as a direct summand, where $V_i$ is the simple representation with $\underline{\dim}V_i=i$. If $i$ is a sink or source, then rep-$\mathcal{S}\langle i\rangle$ is closed under direct summands and extensions. If $i$ is a sink (resp. source), then $\sigma_i^+:\textrm{rep-$\mathcal{S}\langle i\rangle\simeq$ rep-$\sigma_i \mathcal{S}\langle i\rangle$}$ (resp. $\sigma_i^+:\textrm{rep-$\mathcal{S}\langle i\rangle\simeq$ rep-$\sigma_i \mathcal{S}\langle i\rangle$}$) is an equivalence.

\subsection{Construction of Lusztig's symmetries}

Assume $i$ is a sink of $Q$. We first define a map $\mathcal{T}_i$ from $\dot{\mathcal{H}}_q^{\ast}(\Lambda)$ to $\dot{\mathcal{H}}_q^{\ast}(\sigma_i\Lambda)$.

For $\lambda\in\mathcal{P}$, assume that $V_{\lambda}=V_{\lambda_0}\oplus tV_i$ and $V_{\lambda_0}$ contains no direct summand isomorphic to $V_i$. Then Hom$(V_{\lambda_0},V_i)=0$ and Ext$(V_i,V_{\lambda_0})=0$. In this case
\begin{displaymath}
\langle M(\lambda)\rangle=v^{\langle\lambda_0,ti\rangle}u_i^{(t)}\langle M(\lambda_0)\rangle
\end{displaymath}
in $\mathcal{H}_q^{\ast}(\Lambda)$. We define a map  $\mathcal{T}_i:\dot{\mathcal{H}}_q^{\ast}(\Lambda)\rightarrow\dot{\mathcal{H}}_q^{\ast}(\sigma_i\Lambda)$ given by
\begin{equation}\label{Lusztig's symmetries1}
\mathcal{T}_i(\langle M(\lambda)\rangle^+\mathbf{1}_{\zeta}\langle M(\lambda')\rangle^-)
=(-1)^{p_1}v^{q_1}u_{i}^{-(t)}\langle M(\sigma^+_i\lambda_0)\rangle^+\mathbf{1}_{s_i\zeta}u_{i}^{+(t')}\langle M(\sigma^+_i\lambda'_0)\rangle^-
\end{equation}
where $p_1=t+t'-\lambda'_0(h_i)$ and
$q_1=-\langle ti,\lambda_0\rangle-t^2\varepsilon_i+t\varepsilon_i-(\zeta,t\alpha_i)+\langle\lambda'_0,t'i\rangle-(\lambda'_0,i)+t'^2\varepsilon_i-t'\varepsilon_i+(\zeta,t'\alpha_i)$;
\begin{equation}\label{Lusztig's symmetries2}
\mathcal{T}_i(\langle M(\lambda')\rangle^-\mathbf{1}_{\zeta}\langle M(\lambda)\rangle^+)
=(-1)^{p_2}v^{q_2}u_{i}^{+(t')}\langle M(\sigma^+_i\lambda'_0)\rangle^-\mathbf{1}_{s_i\zeta}u_{i}^{-(t)}\langle M(\sigma^+_i\lambda_0)\rangle^+
\end{equation}
where $p_2=t+t'-\lambda'_0(h_i)$ and $q_2=t^2\varepsilon_i+t\varepsilon_i+\langle\lambda_0,ti\rangle-(\zeta,t\alpha_i)-\langle t'i,\lambda'_0\rangle-(\lambda'_0,i)-t'^2\varepsilon_i-t'\varepsilon_i+(\zeta,t'\alpha_i)$.

In fact, the definition of $\mathcal{T}_i$ is induced by the following formulas:
\begin{eqnarray*}
\mathcal{T}_i(\langle M(\lambda)\rangle^+\mathbf{1}_{\zeta})
&=&\langle M(\sigma^+_i\lambda)\rangle^+\mathbf{1}_{s_i\zeta}\\
\mathcal{T}_i(\langle M(\lambda)\rangle^-\mathbf{1}_{\zeta})
&=&(-1)^{\lambda(h_i)}v^{-(\lambda,i)}\langle M(\sigma^+_i\lambda)\rangle^-\mathbf{1}_{s_i\zeta}
\end{eqnarray*}
if $V_{\lambda}$ contains no direct summand isomorphic to $V_i$ and
\begin{eqnarray*}
\mathcal{T}_i(u_i^+\mathbf{1}_{\zeta})
&=&-v^{-(\zeta,\alpha_i)}u_i^-\mathbf{1}_{s_i\zeta}\\
\mathcal{T}_i(u_i^-\mathbf{1}_{\zeta})
&=&-v^{(\zeta,\alpha_i)-2\varepsilon_i}u_i^+\mathbf{1}_{s_i\zeta}.
\end{eqnarray*}

Note that, by the relation (\ref{2}) in the definition of $\dot{\mathcal{H}}^*(\Lambda)$, we can define $\mathcal{T}_i$ on all the generators of $\dot{\mathcal{H}}^*(\Lambda)$. If we can prove that $\mathcal{T}_i$ keeps the relations (\ref{1}) to (\ref{5}), then $\mathcal{T}_i$ induces a map from $\dot{\mathcal{H}}^*(\Lambda)$ to $\dot{\mathcal{H}}^*(\sigma_i\Lambda)$. This is the first main result of this section.

\begin{theorem}\label{Theorem 5D}
Let $i$ be a sink. The formula (\ref{Lusztig's symmetries1}) and (\ref{Lusztig's symmetries2}) induces a $\mathbb{Q}(v)$-algebra isomorphism $\mathcal{T}_i:\dot{\mathcal{H}}^*(\Lambda)\simeq\dot{\mathcal{H}}^*(\sigma_i\Lambda)$
\end{theorem}

The proof of Theorem \ref{Theorem 5D} will be given in the last section.

Let $i$ be a sink. For $j\in I$, if $i=j$, we have $\mathcal{T}_i(u_i^+\mathbf{1}_{\zeta})\in\dot{\mathcal{C}}_q^{\ast}(\sigma_i\Lambda)$ and $\mathcal{T}_i(u_i^-\mathbf{1}_{\zeta})\in\dot{\mathcal{C}}_q^{\ast}(\sigma_i\Lambda)$ since $u_i^+\mathbf{1}_{\zeta}$ and $u_i^-\mathbf{1}_{\zeta}$ are contained in $\dot{\mathcal{C}}_q^{\ast}(\sigma_i\Lambda)$. If $i\neq j$, we have $\mathcal{T}_i(u_j^+\mathbf{1}_{\zeta})=\langle{M(\sigma_i^+(j))}\rangle^+\mathbf{1}_{s_i\zeta}$. Note that $V_{\sigma^+_i(j)}$ is an exceptional object in rep-$\sigma_i\mathcal{S}$. Hence $\langle{M(\sigma_i^+(j))}\rangle\in\dot{\mathcal{C}}_q^{\ast}(\sigma_i\Lambda)$. Hence $\mathcal{T}_i(u_j^+\mathbf{1}_{\zeta})\in\dot{\mathcal{C}}_q^{\ast}(\sigma_i\Lambda)$. Similarly we have $\mathcal{T}_i(u_j^-\mathbf{1}_{\zeta})\in\dot{\mathcal{C}}_q^{\ast}(\sigma_i\Lambda)$. Hence $\mathcal{T}_i$ induces an $\mathbb{Q}(v)$-algebra homomorphism from $\dot{\mathcal{C}}_q^{\ast}(\Lambda)$ to $\dot{\mathcal{C}}_q^{\ast}(\sigma_i\Lambda)$. Note the formula (\ref{Lusztig's symmetries1}) and (\ref{Lusztig's symmetries2}) are independent of the choice of the field. We can consider them as formulas in $\dot{\mathcal{C}}^{\ast}(Q)$ and $\dot{\mathcal{C}}^{\ast}(\sigma_iQ)$. Since both $\dot{\mathcal{C}}^{\ast}(Q)$ and $\dot{\mathcal{C}}^{\ast}(\sigma_iQ)$ are isomorphic to $\dot{U}$, $\mathcal{T}_i$ induces a endomorphism on $\dot{U}$, if we identify $\dot{\mathcal{C}}^{\ast}(Q)$ and $\dot{\mathcal{C}}^{\ast}(\sigma_iQ)$ with $\dot{U}$.

Assume $i$ is a source. For $\lambda\in\mathcal{P}$, assume that $V_{\lambda}=V_{\lambda_0}\oplus tV_i$
and $V_{\lambda_0}$ contains no direct summand isomorphic to $V_i$.
Then Hom$(V_i,V_{\lambda_0})=0$ and Ext$(V_{\lambda_0},V_i)=0$. In this case
\begin{displaymath}
\langle M(\lambda)\rangle=v^{\langle ti,\lambda_0\rangle}\langle M(\lambda_0)\rangle u_i^{(t)}
\end{displaymath}
in $\mathcal{H}_q^{\ast}(\Lambda)$.
We define a map $\mathcal{T}'_i:\dot{\mathcal{H}}_q^{\ast}(\Lambda)\rightarrow\dot{\mathcal{H}}_q^{\ast}(\sigma_i\Lambda)$ given by
\begin{displaymath}
\mathcal{T}'_i(\langle M(\lambda)\rangle^+\mathbf{1}_{\zeta}\langle M(\lambda')\rangle^-)
=(-1)^{p_1}v^{q_1}\langle M(\sigma^+_i\lambda_0)\rangle^+u_{i}^{-(t)}\mathbf{1}_{s_i\zeta}\langle M(\sigma^+_i\lambda'_0)\rangle^-u_{i}^{+(t')}
\end{displaymath}
where $p_1=t-t'-\lambda'_0(h_i)$ and $q_1=\langle ti,\lambda\rangle+t\varepsilon_i+(\zeta,t\alpha_i)-(\lambda'_0,i)-t'\varepsilon_i-t'^2\varepsilon_i-(\zeta,t'\alpha_i)-\langle\lambda'_0,t'i\rangle$;
\begin{displaymath}
\mathcal{T}'_i(\langle M(\lambda')\rangle^-\mathbf{1}_{\zeta}\langle M(\lambda)\rangle^+)
=(-1)^{p_2}v^{q_2}\langle M(\sigma^+_i\lambda'_0)\rangle^-u_{i}^{+(t')}\mathbf{1}_{s_i\zeta}\langle M(\sigma^+_i\lambda_0)\rangle^+u_{i}^{-(t)}
\end{displaymath}
where $p_2=t-t'-\lambda'_0(h_i)$ and $q_2=-t^2\varepsilon_i+t\varepsilon_i+(\zeta,t\alpha_i)-\langle\lambda_0, ti\rangle-(\lambda'_0,i)-t'\varepsilon_i-(\zeta,t'\alpha_i)+\langle t'i,\lambda'\rangle$.

By a similar way, we can prove that $\mathcal{T}'_i$ induces a $\mathbb{Q}(v)$-algebra homomorphism from $\dot{\mathbf{U}}$ to $\dot{\mathbf{U}}$.

Now assume $i$ is a sink of $Q$. Then $i$ is a source of $\sigma_iQ$. We can easily check that $\mathcal{T}_i\mathcal{T}'_i=1$ and $\mathcal{T}'_i\mathcal{T}_i=1$. Hence $\mathcal{T}_i$ is a $\mathbb{Q}(v)$-algebra isomorphism with $\mathcal{T}'_i$ as its inverse.

Hence, we have the following theorem.

\begin{theorem}\label{Theorem 5C}
Let $i$ be a sink. The formula (\ref{Lusztig's symmetries1}) and (\ref{Lusztig's symmetries2}) induces a $\mathbb{Q}(v)$-algebra isomorphism $\mathcal{T}_i: \dot{\mathbf{U}}\simeq\dot{\mathbf{U}}$.
\end{theorem}

Then we will prove that $\mathcal{T}_i$ coincides with $T_i$.

\begin{proposition}[\cite{XY}]\label{proposition5D}
Let $i\neq j\in I$ and $n=a_{ij}$.

(1) If $i$ is a sink, then in $\mathcal{H}_q^{\ast}(\Lambda)$ we have
\begin{displaymath}
\langle M(\lambda)\rangle=\sum^{n}_{t=0}(-1)^tv_i^{-t}u_i^{(t)}u_ju_i^{(n-t)}
\end{displaymath}
where $\lambda\in\mathcal{P}$ is the unique isomorphism class of indecomposable representation with the dimension vector $j+ni$.

(2) If $i$ is a source, then in $\mathcal{H}_q^{\ast}(\Lambda)$ we have
\begin{displaymath}
\langle M(\lambda)\rangle=\sum^{n}_{t=0}(-1)^tv_i^{-t}u_i^{(n-t)}u_ju_i^{(t)}
\end{displaymath}
where $\lambda\in\mathcal{P}$ is the unique isomorphism class of indecomposable representation with the dimension vector $j+ni$.
\end{proposition}

Since $i$ is a sink in $Q$, $i$ is a source in $\sigma_iQ$, and $V_{\sigma_i^+(j)}$ is a unique indecomposable module in rep-$\sigma_i\mathcal{S}$ with dimension vector $j+ni$ where $n=a_{ij}$. Thus by the Proposition \ref{proposition5D},
\begin{displaymath}
\langle M(\sigma_i^+(j))\rangle^+\mathbf{1}_{s_i\zeta}=\sum^{n}_{t=0}(-1)^tv_i^{-t}u_i^{+(n-t)}u_j^+u_i^{+(t)}\mathbf{1}_{s_i\zeta}.
\end{displaymath}
Hence
\begin{displaymath}
\mathcal{T}_i(u_j^+\mathbf{1}_{\zeta})=\sum^{n}_{t=0}(-1)^tv_i^{-t}u_i^{+(n-t)}u_j^+u_i^{+(t)}\mathbf{1}_{s_i\zeta}=T_i(u_j^+\mathbf{1}_{\zeta}).
\end{displaymath}
Similarly we can check $\mathcal{T}_i=T_i$ on other generators.

Hence, we have the following theorem.

\begin{theorem}\label{Theorem 5C}
If $i$ is a sink, then the isomorphism $\mathcal{T}_i:\dot{\mathbf{U}}\rightarrow\dot{\mathbf{U}}$
coincides with $T_i$.
\end{theorem}

\subsection{Braid group relations}\label{11}

Let $A=(a_{ij})_{i,j\in I}$ be a symmetrizable generalized Cartan matrix. If $d(i,j)=a_{ij}a_{ji}\leq 3$, then the order $m(i,j)$ of $s_is_j$ is finite (\cite{Kac}). In fact, we have
$$
m(i,j)=\left\{\begin{array}{c}
         2\,\,\,\,\textrm{if $d(i,j)=0$};\\
         3\,\,\,\,\textrm{if $d(i,j)=1$};\\
         4\,\,\,\,\textrm{if $d(i,j)=2$};\\
         6\,\,\,\,\textrm{if $d(i,j)=3$};\\
         \infty\,\,\,\textrm{if $d(i,j)\geq4$}.
       \end{array}
\right.$$
The braid group of type $A$ is defined by the generators $\{\kappa_i\}_{i\in I}$ and relations
$$\kappa_i\kappa_j\dots=\kappa_j\kappa_i\dots$$
for $i\neq j$ with $m(i,j)\leq+\infty$ factors on both sides, where $m(i,j)$ is the order of $s_is_j$ in $W$, that is,
\begin{eqnarray}
&&\kappa_i\kappa_j=\kappa_j\kappa_i\,\,\,\,\,\,\textrm{if $m(i,j)=2$};\nonumber\\
&&\kappa_i\kappa_j\kappa_i=\kappa_j\kappa_i\kappa_j\,\,\,\,\,\,\textrm{if $m(i,j)=3$};\nonumber\\
&&\kappa_i\kappa_j\kappa_i\kappa_j=\kappa_j\kappa_i\kappa_j\kappa_i\,\,\,\,\,\,\textrm{if $m(i,j)=4$};\nonumber\\
&&\kappa_i\kappa_j\kappa_i\kappa_j\kappa_i\kappa_j\kappa_i=\kappa_j\kappa_i\kappa_j\kappa_i\kappa_j\kappa_i\kappa_j\,\,\,\,\,\,\textrm{if $m(i,j)=6$}.\label{7.1}
\end{eqnarray}

Let $\Lambda$ be a finite dimensional hereditary algebra, and $A$ be the corresponding generalized Cartan matrix. In \cite{XY}, the Lusztig's symmetries on $\mathcal{D}_c(\Lambda)$ are constructed as follows.

\begin{theorem}
Let $i$ be a sink. For all $\lambda\in\mathcal{P}$ and $\mu\in P^{\vee}$, we write $V_{\lambda}\simeq V_{\lambda_0}\oplus tV_i$ where $V_{\lambda_0}$ contain no direct summand isomorphic to $V_i$. Then the map $\tilde{\mathcal{T}}_i$ is defined as follows:
\begin{eqnarray}
&&\tilde{\mathcal{T}}_i(\langle{u_{\lambda}(+)}\rangle)=v^{\langle{\lambda,ti}\rangle}\tilde{K}_{ti}\langle{u_{i}(-)}\rangle^{(t)}\langle{u_{\sigma_i^+\lambda_0}(+)}\rangle;\\
&&\tilde{\mathcal{T}}_i(\langle{u_{\lambda}(-)}\rangle)=v^{\langle{\lambda,ti}\rangle}\tilde{K}_{-ti}\langle{u_{i}(+)}\rangle^{(t)}\langle{u_{\sigma_i^+\lambda_0}(-)}\rangle;\\
&&\tilde{\mathcal{T}}_i(K_{\mu})=K_{s_i(\mu)},
\end{eqnarray}
induces a $\mathbb{Q}(v)$-algebra isomorphism: $\mathcal{D}_c(\Lambda)\simeq\mathcal{D}_c(\sigma_i\Lambda)$.
\end{theorem}

In \cite{XY}, the following theorem is proved.

\begin{theorem}\label{Theorem 7A}
For any $i\neq j\in I$ such that $m=m(i,j)\leq+\infty$, $\tilde{\mathcal{T}}_i$ and $\tilde{\mathcal{T}}_j$ satisfy braid group relations (\ref{7.1}) of type $A$ as maps on $\mathcal{D}_c(\Lambda)$.
\end{theorem}

Let $\Lambda$ be a finite dimensional hereditary algebra. Similarly to the the relation between $\dot{\mathbf{U}}$ and $\mathbf{U}$, We consider the relation between $\dot{\mathcal{H}}_{q}^{\ast}(\Lambda)$ and $\mathcal{D}(\Lambda)$.
For any $\zeta\in P$, we have a surjective linear mapping
\begin{eqnarray*}
\pi_{\zeta}:\mathcal{D}(\Lambda)&\rightarrow&\dot{\mathcal{H}}_{q}^{\ast}(\Lambda)\mathbf{1}_{\zeta}\\
\langle{u_{\alpha}(+)}\rangle\langle{u_{\beta}(-)}\rangle K_{\mu}&\mapsto&(-1)^{tr(\beta)}v^{m(\beta)+\zeta(\mu)}\langle{M(\alpha)}\rangle^+\langle{M(\beta)}^-\rangle\mathbf{1}_{\zeta}
\end{eqnarray*}
where $\beta=\sum_{i\in I}b_i\alpha_i$, $tr(\beta)=\sum_{i\in I}b_i$ and $m(\beta)=\sum_{i\in I}b_i\varepsilon_i$.
The kernel of $\pi_{\zeta}$ is $$\sum_{\mu\in P^{\vee}}\mathcal{D}(\Lambda)(K_{\mu}-v^{\zeta(\mu)}).$$ For any $\zeta,\zeta'\in P$, $\beta\in\mathbb{Z}I$ and any $x\in \mathcal{D}(\Lambda),y\in \mathcal{D}(\Lambda)(\beta)$,
\begin{displaymath}
\pi_{\zeta}(x)\pi_{\zeta'}(y)=
\left\{
  \begin{array}{cc}
    \pi_{\zeta'}(xy) & \textrm{if}\,\,\,\zeta=\zeta'+\beta\\
     0 & \textrm{otherwise}
  \end{array}
\right..
\end{displaymath}

Our main result in this subsection is the following.

\begin{theorem}
Let $\Lambda$ be a finite dimensional hereditary algebra, and $A$ be the corresponding generalized Cartan matrix. For any $i\neq j\in I$ such that $m=m(i,j)\leq+\infty$, $\mathcal{T}_i$ and $\mathcal{T}_j$ satisfy braid group relations (\ref{7.1}) of type $A$ as maps on $\dot{\mathcal{C}}_q^{\ast}(\Lambda)$.
\end{theorem}
\begin{prof}
For all $\lambda\in\mathcal{P}$ and $\mu\in P^{\vee}$, we write $V_{\lambda}\simeq V_{\lambda_0}\oplus tV_i$ where $V_{\lambda_0}$ contain no direct summand isomorphic to $V_i$. We need to check that for any $\zeta\in P$
\begin{eqnarray}
&&\pi_{s_i\zeta}(\tilde{\mathcal{T}}_i(\langle{u_{\lambda}(+)}\rangle))=\mathcal{T}_i(\pi_{\zeta}(\langle{u_{\lambda}(+)}\rangle));\label{a1}\\
&&\pi_{s_i\zeta}(\tilde{\mathcal{T}}_i(\langle{u_{\lambda}(-)}\rangle))=\mathcal{T}_i(\pi_{\zeta}(\langle{u_{\lambda}(-)}\rangle));\label{a2}\\
&&\pi_{s_i\zeta}(\tilde{\mathcal{T}}_i(K_{\mu}))=\mathcal{T}_i(\pi_{\zeta}(K_{\mu})).\label{a3}
\end{eqnarray}
First
\begin{eqnarray*}
\pi_{s_i\zeta}(\tilde{\mathcal{T}}_i(\langle{u_{\lambda}(+)}\rangle))
&=&\pi_{s_i\zeta}(v^{\langle{\lambda,ti}\rangle}\tilde{K}_{ti}\langle{u_{i}(-)}\rangle^{(t)}\langle{u_{\sigma_i^+\lambda_0}(+)}\rangle)\\
&=&v^{\langle{\lambda,ti}\rangle+(\sigma_i^+\lambda_0-ti,ti)}
\pi_{s_i\zeta}(\langle{u_{i}(-)}\rangle^{(t)}\langle{u_{\sigma_i^+\lambda_0}(+)}\rangle\tilde{K}_{ti})\\
&=&v^{\langle{\lambda,ti}\rangle+(\sigma_i^+\lambda_0-ti,ti)+(s_i\zeta,t\alpha_i)}(-1)^{t}v^{m(ti)}
u_{i}^{-(t)}\langle{M(\sigma_i^+\lambda_0)}\rangle^+\mathbf{1}_{s_i\zeta}\\
&=&(-1)^{t}v^{-\langle{ti,\lambda_0}\rangle-t^2\varepsilon_i+t\varepsilon_i-(\zeta,t\alpha_i)}
u_{i}^{-(t)}\langle{M(\sigma_i^+\lambda_0)}\rangle^+\mathbf{1}_{s_i\zeta}\\
&=&\mathcal{T}_i(\langle{M(\sigma_i^+\lambda_0)}\rangle^+\mathbf{1}_{\zeta})\\
&=&\mathcal{T}_i(\pi_{\zeta}(\langle{u_{\lambda}(+)}\rangle)).
\end{eqnarray*}
Hence we have formula (\ref{a1}). Similarly, we can get formula (\ref{a2}) and (\ref{a3}). Then Theorem \ref{Theorem 7A} implies this theorem.
\end{prof}

\section{Lusztig's symmetries on the modified form of Ringel-Hall algebras}

\subsection{The structure of Ringel-Hall algebras}\label{12}

First we recall the structure of the Ringel-Hall algebra considered in \cite{SV2} and \cite{DX On double Ringel¨CHall algebras}.

We consider a bilinear form $\psi:\mathcal{H}_q^{\ast}(\Lambda)\times\mathcal{H}_q^{\ast}(\Lambda)$ as
\begin{displaymath}
\psi(\langle{M(\beta)}\rangle,\langle{M(\beta')}\rangle)=\frac{|V_{\beta}|}{a_{\beta}}\delta_{\beta\beta'}
\end{displaymath}
for $\beta,\beta'\in\mathcal{P}$.

Let $\mathfrak{d}_0(\Lambda)=\mathcal{C}_q^{\ast}(\Lambda)$. We can define $\mathfrak{d}_m(\Lambda)$ and $L_{\pi_m}(\Lambda)$ inductively. For $m\geq1$, assume $\mathfrak{d}_{m-1}(\Lambda)$ has been constructed. Let $\pi_m\in\mathbb{Z}I$ have smallest trace such that $\mathfrak{d}_{m-1}(\Lambda)_{\pi_m}\neq\mathcal{H}_q^{\ast}(\Lambda)_{\pi_m}$. Then $L_{\pi_m}(\Lambda)$ is defined as follow:
\begin{displaymath}
L_{\pi_m}(\Lambda):=\{x\in\mathcal{H}_q^{\ast}(\Lambda)_{\pi_m}|\psi(x,\mathfrak{d}_{m-1}(\Lambda)_{\pi_m})=0\}.
\end{displaymath}
We define $\mathfrak{d}_m(\Lambda)$ as the subalgebra of $\mathcal{H}_q^{\ast}(\Lambda)$ generated by $\mathfrak{d}_{m-1}(\Lambda)$ and $L_{\pi_m}(\Lambda)$. Hence there is a chain of subalgebras of $\mathcal{H}_q^{\ast}(\Lambda)$
\begin{displaymath}
\mathfrak{d}_{0}(\Lambda)\subset\mathfrak{d}_{1}(\Lambda)\subset\ldots\mathfrak{d}_{m}(\Lambda)\subset\ldots\subset\mathcal{H}_q^{\ast}(\Lambda).
\end{displaymath}

For $m\geq1$, let $\eta_m=\dim{L_{\pi_m}}$. There exists a bases $\{x_{(m,p)}|1\leq p\leq\eta_m\}$ of $L_{\pi_m}$ and nonzero numbers $\chi_{(m,p)}\in\mathbb{Q}(v),1\leq p\leq\eta_m$ such that
\begin{displaymath}
\psi(x_{(m,p)},\chi_{(m,p)}x_{(m,q)})=\frac{-1}{v-v^{-1}}\delta_{pq}.
\end{displaymath}
Set $x_i=u_i$ and $J=\{(m,p)|m\geq 1,1\leq p\leq\eta_m\}$. The elements in the set $\{x_j|j\in I\cup J\}$ generate the Ringel-Hall algebra $\mathcal{H}_q^{\ast}(\Lambda)$.

Let $y_i=-v_i^{-1}u_i$ for all $i\in I$ and $y_j=\chi_jx_j$ for all $j\in J$. By \cite{SV2} and \cite{DX On double Ringel¨CHall algebras}, the double Ringel-Hall algebra $\mathcal{D}(\Lambda)$ is generated by the elements $x_i(+),y_i(-),i\in I\cup J$ and $K_{\mu},\mu\in P^{\vee}$ subject to the following relations:
\begin{equation}
K_{0}=\mathbf{1}, K_{\mu}K_{\mu'}=K_{\mu+\mu'}\,\,\,\textrm{for all $\mu,\mu'\in P^{\vee}$};
\end{equation}
\begin{equation}
K_{\mu}x_{i}(+)K_{-\mu}=v^{\delta_i(\mu)}x_i(+)\,\,\,\textrm{for all $i\in I\cup J$, $\mu\in P^{\vee}$};
\end{equation}
\begin{equation}
K_{\mu}y_{i}(-)K_{-\mu}=v^{-\delta_i(\mu)}y_i(-)\,\,\,\textrm{for all $i\in I\cup J$, $\mu\in P^{\vee}$};
\end{equation}
\begin{equation}
x_i(+)y_j(-)-y_j(-)x_i(+)=\delta_{ij}\frac{{\tilde{K}}_{\delta_i}-{\tilde{K}}_{-\delta_i}}{v_i-v_i^{-1}}\,\,\,\textrm{for all $i,j\in I\cup J$};
\end{equation}
for $i\in I$, $j\in I\cup J$ and $i\neq j$, setting $b=1-a_{ij}$,
\begin{equation}
\sum_{k=0}^{b}(-1)^{k}x_i(+)^{(k)}x_j(+)x_i(+)^{(b-k)}=0,
\end{equation}
and
\begin{equation}
\sum_{k=0}^{b}(-1)^{k}y_i(-)^{(k)}y_j(-)y_i(-)^{(b-k)}=0;
\end{equation}
for any $i,j\in I\cup J$ with $(\delta_i,\delta_j)=0$,
\begin{equation}
x_i(+)x_j(+)=x_j(+)x_i(+),\,\,\,\,\,\,\,\,y_i(-)y_j(-)=y_j(-)y_i(-).
\end{equation}
Here, $\delta_i=\alpha_i$ for $i\in I$, $\delta_j=\pi_m$ for $j=(m,p)\in J$ and $a_{ij}=2\frac{(\delta_i,\delta_j)}{(\delta_i,\delta_i)}$.

Note that $\tilde{A}=(a_{ij})_{ij\in I\cup J}$ is a Borcherds-Cartan matrix. We can define a modified quantized enveloping algebra $\dot{\mathbf{U}}(\tilde{A})$ of the generalized Kac-Moody algebra associated to $\tilde{A}$. $\dot{\mathbf{U}}(\tilde{A})$ is generated by the elements $E_i\mathbf{1}_{\zeta},F_i\mathbf{1}_{\zeta}$ for all $i\in I\cup J$ and $\zeta\in P$ subject to the following relations:
\begin{equation}
\mathbf{1}_{\zeta}\mathbf{1}_{\zeta'}=\delta_{\zeta\zeta'}\mathbf{1}_{\zeta}\,\,\,\textrm{for all $\zeta,\zeta'\in P$};\label{31}
\end{equation}
\begin{equation}
E_i\mathbf{1}_{\zeta}=\mathbf{1}_{\zeta+\delta_i}E_i,F_i\mathbf{1}_{\zeta}=\mathbf{1}_{\zeta-\delta_i}F_i\,\,\,\textrm{for all $i\in I\cup J$, $\zeta\in P$};
\end{equation}
\begin{equation}
(E_i\mathbf{1}_{\zeta-\delta_j})(F_j\mathbf{1}_{\zeta})
-(F_j\mathbf{1}_{\zeta+\delta_i})(E_i\mathbf{1}_{\zeta})
=\delta_{ij}(-1)^{tr\delta_j}v^{-m(\delta_j)}\frac{v^{(\zeta,\delta_i)}-v^{-(\zeta,\delta_i)}}{v_i-v_i^{-1}}\,\,\,\textrm{for all $i,j\in I\cup J$};
\end{equation}
for $i\in I$, $j\in I\cup J$ and $i\neq j$, setting $b=1-a_{ij}$,
\begin{equation}
\sum_{k=0}^{b}(-1)^{k}(E_i^{(k)}\mathbf{1}_{\zeta+(b-k)\delta_i+\delta_j})
(E_j\mathbf{1}_{\zeta+(b-k)\delta_i})(E_i^{(b-k)}\mathbf{1}_{\zeta})=0,
\end{equation}
and
\begin{equation}
\sum_{k=0}^{b}(-1)^{k}(F_i^{(k)}\mathbf{1}_{\zeta-(b-k)\delta_i-\delta_j})
(F_j\mathbf{1}_{\zeta-(b-k)\delta_i})(F_i^{(b-k)}\mathbf{1}_{\zeta})=0;
\end{equation}
for any $i,j\in I\cup J$ with $(\delta_i,\delta_j)=0$,
\begin{equation}
(E_i\mathbf{1}_{\zeta+\delta_j})(E_j\mathbf{1}_{\zeta})=(E_j\mathbf{1}_{\zeta+\delta_i})(E_i\mathbf{1}_{\zeta})
,\,\,\,\,\,\,\,\,(F_i\mathbf{1}_{\zeta-\delta_j})(F_j\mathbf{1}_{\zeta})=(F_j\mathbf{1}_{\zeta-\delta_i})(F_i\mathbf{1}_{\zeta}),\label{32}
\end{equation}
where
$$E_i^{(k)}\mathbf{1}_{\zeta}=\frac{1}{[k]_{v_i}!}\prod_{s=1}^{k}E_i\mathbf{1}_{\zeta+(k-s)\delta_i},$$
$$F_i^{(k)}\mathbf{1}_{\zeta}=\frac{1}{[k]_{v_i}!}\prod_{s=1}^{k}F_i\mathbf{1}_{\zeta-(k-s)\delta_i}.$$

Since there exists a map $\pi_{\zeta}:\mathcal{D}(\Lambda)\rightarrow\dot{\mathcal{H}}_{q}^{\ast}(\Lambda)\mathbf{1}_{\zeta}$ for any $\zeta\in P$, the algebra $\dot{\mathcal{H}}_q^{\ast}(\Lambda)$ is generated by the elements $x_i^+\mathbf{1}_{\zeta},y_i^-\mathbf{1}_{\zeta}$ for all $i\in I\cup J$ and $\zeta\in P$ subject to the relations (\ref{31}) to (\ref{32}). Hence, we have an isomorphism $\iota:\dot{\mathcal{H}}_{q}^{\ast}(\Lambda)\simeq \dot{\mathbf{U}}(\tilde{A})$ mapping $x_i^+\mathbf{1}_{\zeta}$ (resp. $y_i^-\mathbf{1}_{\zeta}$) to $E_i\mathbf{1}_{\zeta}$ (resp. $F_i\mathbf{1}_{\zeta}$).


There is an operator $\tau$ on $\mathcal{H}_q^{\ast}(\Lambda)$ defined as follows:
\begin{eqnarray*}
\tau{\langle{M(\lambda)}\rangle}&=&(-1)^{tr\alpha}v^{-\tau(\alpha)}\\
&\times&\left(\delta_{\lambda0}
+\sum_{m\geq1}(-1)^{m}\sum_{\pi\in\mathcal{P},\lambda_1,\ldots,\lambda_m\in\mathcal{P}\backslash\{0\}}
v^{2\sum_{i<j}\langle{\lambda_i,\lambda_j}\rangle}\times\right.\\
&&\left.\frac{a_{\lambda_1\ldots a_{\lambda_m}}}{a_{\lambda}}g_{\lambda_1,\ldots,\lambda_m}^{\lambda}g_{\pi}^{\lambda_1,\ldots,\lambda_m}\langle{M(\pi)}\rangle\right)
\end{eqnarray*}
where $\lambda\in\mathcal{P}$, $u_{\lambda}\in\mathcal{H}_q^{\ast}(\Lambda)_{\alpha},\alpha=\sum_{i}k_i\alpha_i\in\mathbb{N}[I]$,
$tr\alpha=\sum_{i}k_i$ and $\tau(\alpha)=((\alpha,\alpha)-\sum_ik_i(i,i))/2$.

\subsection{Lusztig's symmetries on the modified form of the Ringel-Hall algebras}

We first recall the definition of Lusztig's symmetries of $\mathcal{D}(\Lambda)$ defined in \cite{DX On double Ringel¨CHall algebras}. For all $i\in I$, define $\tilde{T}_i:\mathcal{D}(\Lambda)\rightarrow\mathcal{D}(\Lambda)$ on generators as follows
\begin{eqnarray*}
&&\tilde{T}_i(x_i(+))=-y_i(-)\tilde{K}_{i},\tilde{T}_i(y_i(-))=-\tilde{K}_{-i}x_i(+);\\
&&\tilde{T}_i(x_j(+))=\sum_{r+s=-a_{ij}}(-1)^rv_i^{-r}x_i(+)^{(s)}x_j(+)x_i(+)^{(r)} \textrm{for $i\neq j\in I\cup J$};\\
&&\tilde{T}_i(y_j(-))=\sum_{r+s=-a_{ij}}(-1)^rv_i^{r}y_i(-)^{(r)}y_j(-)y_i(-)^{(s)} \textrm{for $i\neq j\in I\cup J$};\\
&&\tilde{T}_i(K_{\mu})=K_{\mu-\alpha_{i}(\mu)h_i} \textrm{for $\mu\in P^{\vee}$}.
\end{eqnarray*}

Under the maps
$$\pi_{\zeta}:\mathcal{D}(\Lambda)\rightarrow\dot{\mathcal{H}}_{q}^{\ast}(\Lambda)\mathbf{1}_{\zeta},$$
Lusztig's symmetries $\tilde{T}_i$ of $\mathcal{D}(\Lambda)$ induce Lusztig's symmetries $T_i:\dot{\mathcal{H}}_{q}^{\ast}(\Lambda)\rightarrow\dot{\mathcal{H}}_{q}^{\ast}(\Lambda)$. From the formulas above, we get
\begin{eqnarray*}
&&T_i(x_i^+\mathbf{1}_{\zeta})=-v_i^{-\zeta(h_i)}\tilde{y}_i^-\mathbf{1}_{s_i\zeta}\,\,\,\textrm{for $\zeta\in P$};\\
&&T_i(\tilde{y}_i^-\mathbf{1}_{\zeta})=-v_i^{-(2-\zeta(h_i))}x_i^+\mathbf{1}_{s_i\zeta}\,\,\,\textrm{for $\zeta\in P$};\\
&&T_i(x_j^+\mathbf{1}_{\zeta})=\sum_{r+s=-a_{ij}}(-1)^rv_i^{-r}x_i^{+(s)}x_j^+x_i^{+(r)}\mathbf{1}_{s_i\zeta}\,\,\,\textrm{for $i\neq j\in I\cup J$};\\
&&T_i(\tilde{y}_j^-\mathbf{1}_{\zeta})=\sum_{r+s=-a_{ij}}(-1)^rv_i^{r}\tilde{y}_i^{-(r)}\tilde{y}_j^-\tilde{y}_i^{-(s)}\mathbf{1}_{s_i\zeta}\,\,\,\textrm{for $i\neq j\in I\cup J$}
\end{eqnarray*}
where $\tilde{y}_i=(-1)^{tr\delta_i}v^{m(\delta_i)}y_i$ for all $i\in I\cup J$. Note that
$\pi_{\zeta}(y_i(-))=\tilde{y}_i^-\mathbf{1}_{\zeta}$.


We define
\begin{displaymath}
\psi^{\pm}_{\zeta}(x^{\pm}\mathbf{1}_{\zeta},x'^{\pm}\mathbf{1}_{\zeta})=\psi(x,x')
\end{displaymath}
for every $\zeta\in P$. Let $\mathcal{H}_q^{\ast}(\Lambda)\langle{i}\rangle$ be the subspace of $\mathcal{H}_q^{\ast}(\Lambda)$ spanned by the elements in the set
$$\{\langle{M(\alpha)}\rangle|\alpha\in\mathcal{P},V_{\alpha}\in\textrm{rep-}\mathcal{S}\langle{i}\rangle\}$$
and $\mathfrak{d}_{m}(\Lambda)\langle{i}\rangle=\mathfrak{d}_{m}(\Lambda)\cap\mathcal{H}_q^{\ast}(\Lambda)\langle{i}\rangle$.

\begin{proposition}\label{proposition 6A}
Let $i\in I$ be a sink. For all $\mu\in P$ and all $x,x'\in\mathcal{H}_q^{\ast}(\Lambda)\langle{i}\rangle_{\mu}$, we have
\begin{displaymath}
\psi^{\pm}_{\zeta}(x^{\pm}\mathbf{1}_{\zeta},x'^{\pm}\mathbf{1}_{\zeta})
=\psi^{\pm}_{s_i\zeta}(T_i(x^{\pm}\mathbf{1}_{\zeta}),T_i(x'^{\pm}\mathbf{1}_{\zeta}))
\end{displaymath}
\end{proposition}
\begin{prof}
In \cite{DX Ringel¨CHall algebras and Lusztig's symmetries}, it is proved that
$$
\psi(x,x')=\psi(\tilde{T}_i(x),\tilde{T}_i(x')).
$$
From the definition of $\psi^{\pm}_{\zeta}(-,-)$,
\begin{eqnarray*}
&&\psi^{\pm}_{s_i\zeta}(T_i(x^{\pm}\mathbf{1}_{\zeta}),T_i(x'^{\pm}\mathbf{1}_{\zeta}))\\
&=&\psi^{\pm}_{s_i\zeta}(\tilde{T}_i(x)^{\pm}\mathbf{1}_{s_i\zeta},\tilde{T}_i(x')^{\pm}\mathbf{1}_{s_i\zeta})\\
&=&\psi(\tilde{T}_i(x),\tilde{T}_i(x'))\\
&=&\psi(x,x')\\
&=&\psi^{\pm}_{\zeta}(x^{\pm}\mathbf{1}_{\zeta},x'^{\pm}\mathbf{1}_{\zeta}).
\end{eqnarray*}
\end{prof}


\subsection{Relation between the Lusztig's symmetries and the BGP-reflection functors}

In this subsection, we consider the relation between the Lusztig's symmetries and the BGP-reflection functors. The method is similar to these in \cite{DX Ringel¨CHall algebras and Lusztig's symmetries}.

\begin{proposition}\label{proposition 6B}
Let $i\in I$ be a sink. For each $x,x'\in\mathcal{H}_q^{\ast}(\Lambda)\langle{i}\rangle$, we have
\begin{displaymath}
\psi^{\pm}_{\zeta}(x^{\pm}\mathbf{1}_{\zeta},x'^{\pm}\mathbf{1}_{\zeta})=
\psi^{\pm}_{s_i\zeta}(\mathcal{T}_i(x^{\pm}\mathbf{1}_{s_i\zeta}),\mathcal{T}_i(x'^{\pm}\mathbf{1}_{s_i\zeta})).
\end{displaymath}
\end{proposition}
\begin{prof}
Let $V_{\beta},V_{\beta'}\in\textrm{rep-}Q\langle{i}\rangle$. Then
\begin{eqnarray*}
&&\psi^{+}_{s_i\zeta}(\mathcal{T}_i(\langle{M(\beta)}\rangle^{+}\mathbf{1}_{\zeta}),\mathcal{T}_i(\langle{M(\beta')}\rangle^{+}\mathbf{1}_{\zeta}))\\
&=&\psi^{+}_{s_i\zeta}
(\langle{M(\sigma_i^+\beta)}\rangle^{+}\mathbf{1}_{s_i\zeta},\langle{M(\sigma_i^+\beta')}\rangle^{+}\mathbf{1}_{s_i\zeta})\\
&=&\psi(\langle{M(\sigma_i^+\beta)}\rangle,\langle{M(\sigma_i^+\beta')})\\
&=&\frac{|V_{\sigma_i^+\beta}|}{a_{\sigma_i^+\beta}}\delta_{\sigma_i^+\beta\sigma_i^+\beta'}\\
&=&\frac{|V_{\beta}|}{a_{\beta}}\delta_{\beta\beta'}\\
&=&\psi(\langle{M(\beta)}\rangle,\langle{M(\beta')})\\
&=&\psi^{+}_{\zeta}
(\langle{M(\beta)}\rangle^{+}\mathbf{1}_{\zeta},\langle{M(\beta')}\rangle^{+}\mathbf{1}_{\zeta}).
\end{eqnarray*}
Hence we have
\begin{displaymath}
\psi^{+}_{\zeta}(x^{+}\mathbf{1}_{\zeta},x'^{+}\mathbf{1}_{\zeta})
=\psi^{+}_{s_i\zeta}(\mathcal{T}_i(x^{+}\mathbf{1}_{s_i\zeta}),\mathcal{T}_i(x'^{+}\mathbf{1}_{s_i\zeta})).
\end{displaymath}

Similarly we can prove that
\begin{displaymath}
\psi^{-}_{\zeta}(x^{-}\mathbf{1}_{\zeta},x'^{-}\mathbf{1}_{\zeta})
=\psi^{-}_{s_i\zeta}(\mathcal{T}_i(x^{-}\mathbf{1}_{s_i\zeta}),\mathcal{T}_i(x'^{-}\mathbf{1}_{s_i\zeta})).
\end{displaymath}
\end{prof}

\begin{theorem}
Let $i\in I$ be a sink. Then for each $m\geq 1$, $\mathcal{T}_iT_i^{-1}$ induces bijective maps from $L_{\pi_m}(\Lambda)^{\pm}\mathbf{1}_{\zeta}$ to $ L_{\pi_m}(\sigma_i\Lambda)^{\pm}\mathbf{1}_{\zeta}$.
\end{theorem}
\begin{prof}
We first prove the theorem for $L^{+}_{\pi_m}(\Lambda)\mathbf{1}_{\zeta}$.
By the definition we have
\begin{displaymath}
L_{\pi_m}(\Lambda)=\{x\in\mathcal{H}_q^{\ast}(\Lambda)_{\pi_m}|\psi(x,\mathfrak{d}_{m-1}(\Lambda)_{\pi_m})=0\}.
\end{displaymath}
By \cite{DX Ringel¨CHall algebras and Lusztig's symmetries}, we have $L_{\pi_m}(\Lambda)\subset{}^{\tau}\mathcal{H}_q^{\ast}(\Lambda)\langle{i}\rangle_{\pi_m}$,  $\mathfrak{d}_{m-1}(\Lambda)\langle{i}\rangle=\sum_{s\geq1}{}^{\tau}\mathfrak{d}_{m-1}(\Lambda)\langle{i}\rangle x^s_i$ and $\psi(x,{}^{\tau}\mathfrak{d}_{m-1}(\Lambda)\langle{i}\rangle x^s_i)=0$ for $x\in{}^{\tau}\mathcal{H}_q^{\ast}(\Lambda)\langle{i}\rangle$, where
${}^{\tau}\mathcal{H}_q^{\ast}(\Lambda)\langle{i}\rangle:=\tau(\mathcal{H}_q^{\ast}(\Lambda)\langle{i}\rangle)$ and $^{\tau}\mathfrak{d}_{m}(\Lambda)\langle{i}\rangle:=\tau(\mathfrak{d}_{m}(\Lambda)\langle{i}\rangle)$.
Then we have
\begin{displaymath}
L_{\pi_m}(\Lambda)=\{x\in{}^{\tau}\mathcal{H}_q^{\ast}(\Lambda)\langle{i}\rangle_{\pi_m}|\psi(x,{}^{\tau}\mathfrak{d}_{m-1}(\Lambda)\langle{i}\rangle_{\pi_m})=0\}.
\end{displaymath}

We have the following isomorphisms
\begin{displaymath}
{}^{\tau}\mathfrak{d}_{m-1}(\Lambda)\langle{i}\rangle_{\pi_m}^{+}\mathbf{1}_{\zeta}
\xrightarrow{T_i^{-1}}\mathfrak{d}_{m-1}(\Lambda)\langle{i}\rangle_{s_i\pi_m}^{+}\mathbf{1}_{s_i\zeta}
\xrightarrow{\mathcal{T}_i}\mathfrak{d}_{m-1}(\sigma_i^+\Lambda)\langle{i}\rangle_{\pi_m}^{+}\mathbf{1}_{\zeta}.
\end{displaymath}
The first isomorphism is showed in \cite{DX On double Ringel¨CHall algebras}. For the second one, we have proved that $\mathcal{T}_i$ is an isomorphism in Theorem \ref{Theorem 5D}. Hence we just need to show
\begin{displaymath}
\mathcal{T}_i(\mathfrak{d}_{m}(\Lambda)\langle{i}\rangle_{\pi_m}^{+}\mathbf{1}_{\zeta})\subset \mathfrak{d}_{m}(\sigma_i^+\Lambda)\langle{i}\rangle_{s_i\pi_m}^+\mathbf{1}_{s_i\zeta}.
\end{displaymath}
By \cite{DX On double Ringel¨CHall algebras}, we know
\begin{displaymath} \mathfrak{d}_{m}(\sigma_i^+\Lambda)\langle{i}\rangle_{s_i\pi_m}=\mathcal{H}^{\ast}(\sigma_i^+\Lambda)\langle{i}\rangle_{s_i\pi_m}.
\end{displaymath}
Hence we have
\begin{displaymath}
\mathcal{T}_i(\mathfrak{d}_{m}(\Lambda)\langle{i}\rangle_{\pi_m}^{+}\mathbf{1}_{\zeta})\subset\mathcal{H}^{\ast}(\sigma_i^+\Lambda)\langle{i}\rangle_{s_i\pi_m}^+\mathbf{1}_{s_i\zeta}
=\mathfrak{d}_{m}(\sigma_i^+\Lambda)\langle{i}\rangle_{s_i\pi_m}^+\mathbf{1}_{s_i\zeta}.
\end{displaymath}

Take any $x\in L_{\pi_m}(\Lambda)$. Then $\psi(x,{}^{\tau}\mathfrak{d}_{m-1}(\Lambda)\langle{i}\rangle_{\pi_m})=0$. By Proposition \ref{proposition 6A} and Proposition \ref{proposition 6B} we have
\begin{eqnarray*}
0&=&\psi(x,{}^{\tau}\mathfrak{d}_{m-1}(\Lambda)\langle{i}\rangle_{\pi_m})\\
&=&\psi^+_{\lambda}(x^+\mathbf{1}_{\zeta},{}^{\tau}\mathfrak{d}_{m-1}(\Lambda)\langle{i}\rangle^+_{\pi_m}\mathbf{1}_{\zeta})\\
&=&\psi^+_{\lambda}(\mathcal{T}_iT_i^{-1}(x^+\mathbf{1}_{\zeta}),\mathcal{T}_iT_i^{-1}({}^{\tau}\mathfrak{d}_{m-1}(\Lambda)\langle{i}\rangle^+_{\pi_m}\mathbf{1}_{\zeta}))\\
&=&\psi^+_{\lambda}(\mathcal{T}_iT_i^{-1}(x^+\mathbf{1}_{\zeta}),{}^{\tau}\mathfrak{d}_{m-1}(\sigma_i^+\Lambda)\langle{i}\rangle^+_{\pi_m}\mathbf{1}_{\zeta}).
\end{eqnarray*}
Hence $\mathcal{T}_iT_i^{-1}(x^+\mathbf{1}_{\zeta})\in L_{\pi_m}(\Lambda)^+\mathbf{1}_{\zeta}$. Conversely, $\mathcal{T}_iT_i^{-1}(x^+\mathbf{1}_{\zeta})\in L_{\pi_m}(\Lambda)^+\mathbf{1}_{\zeta}$ implies $x^+\mathbf{1}_{\zeta}\in L_{\pi_m}(\Lambda)^+\mathbf{1}_{\zeta}$. Hence $\mathcal{T}_iT_i^{-1}$ induces bijective maps from $L^{+}_{\pi_m}(\Lambda)\mathbf{1}_{\zeta}$ to $L^{+}_{\pi_m}(\sigma_i\Lambda)\mathbf{1}_{\zeta}$.

Similarly, we can prove that $\mathcal{T}_iT_i^{-1}$ induces bijective maps from $L^{-}_{\pi_m}(\Lambda)\mathbf{1}_{\zeta}$ to
$L^{-}_{\pi_m}(\sigma_i\Lambda)\mathbf{1}_{\zeta}$.
\end{prof}

As in Section \ref{12}, by choosing the basis $\{x_{(m,p)}|1\leq p\leq\eta_m\}$ of $L_{\pi_m}$ for all $m$, we get a set of generators $G=\{x_i^+\mathbf{1}_{\zeta},y_i^-\mathbf{1}_{\zeta}|i\in I\cup J,\zeta\in P^{\vee}\}$ of $\dot{\mathcal{H}}^{\ast}_q(\Lambda)$ and $\dot{\mathcal{H}}^{\ast}_q(\Lambda)$ is generated by these elements subject to the relations (\ref{31}) to (\ref{32}). If $i\in I$ is a sink, the theorem above implies that the image of $G$ under $\mathcal{T}_iT_i^{-1}$ becomes a set of generators of $\dot{\mathcal{H}}^{\ast}_q(\sigma_i\Lambda)$ subject to the same relations. Hence, we also have an isomorphism $\iota':\dot{\mathcal{H}}_{q}^{\ast}(\sigma_i\Lambda)\simeq \dot{\mathbf{U}}(\tilde{A})$ mapping $\mathcal{T}_iT_i^{-1}(x_i^+\mathbf{1}_{\zeta})$ (resp. $\mathcal{T}_iT_i^{-1}(y_i^-\mathbf{1}_{\zeta})$) to $E_i\mathbf{1}_{\zeta}$ (resp. $F_i\mathbf{1}_{\zeta}$). Under the isomorphisms $\iota$ and $\iota'$, the maps $\mathcal{T}_i$ and $T_i$ induce maps on $\dot{\mathbf{U}}(\tilde{A})$, which are also denoted by $\mathcal{T}_i$ and $T_i$ respectively. Then we have the following theorem.

\begin{theorem}
Let $i\in I$ be a sink. Then the isomorphisms $\mathcal{T}_i$ and $T_i$ coincide as maps from $\dot{\mathbf{U}}(\tilde{A})$ to $\dot{\mathbf{U}}(\tilde{A})$.
\end{theorem}
\begin{prof}
Under the isomorphisms $\iota$ and $\iota'$, we get a map $\mathcal{T}_iT_i^{-1}$ from $\dot{\mathbf{U}}(\tilde{A})$ to $\dot{\mathbf{U}}(\tilde{A})$. Note that $\mathcal{T}_iT_i^{-1}$ sends the generators $E_i\mathbf{1}_{\zeta}$ and $F_i\mathbf{1}_{\zeta}$ to themselves. Hence $\mathcal{T}_iT_i^{-1}$ is the identical map on $\dot{\mathbf{U}}(\tilde{A})$. So $\mathcal{T}_i$ and $T_i$ coincide.
\end{prof}

\subsection{Braid group relations}

In \cite{DX Ringel¨CHall algebras and Lusztig's symmetries}, the following theorem is proved.

\begin{theorem}
For any $i\neq j\in I$ such that $m=m(i,j)\leq+\infty$, $\tilde{\mathcal{T}}_i$ and $\tilde{\mathcal{T}}_j$ satisfy braid group relations (\ref{7.1}) of type $A$ as maps on $\mathcal{D}(\Lambda)$.
\end{theorem}

Similarly to the case in Section \ref{11}, we have the following theorem.

\begin{theorem}
Let $\Lambda$ be a finite dimensional hereditary algebra, and $A$ be the corresponding generalized Cartan matrix. For any $i\neq j\in I$ such that $m=m(i,j)\leq+\infty$, $\mathcal{T}_i$ and $\mathcal{T}_j$ satisfy braid group relations (\ref{7.1}) of type $A$ as maps on $\dot{\mathbf{U}}(\tilde{A})$.
\end{theorem}

\section{The proof of Theorem \ref{Theorem 5D}}

Let $i$ be a sink and we follow the method used in \cite{XY}.

From the definition of $\mathcal{T}_i$, we have the following proposition.
\begin{proposition}\label{proposition5B}
For any $\lambda,\lambda'\in\mathcal{P}$, we have
\begin{equation}
\mathcal{T}_i(\langle{M(\lambda)}\rangle^+\mathbf{1}_{\zeta})=\mathcal{T}_i(\mathbf{1}_{\zeta+\lambda}\langle{M(\lambda)}\rangle^+),
\mathcal{T}_i(\langle{M(\lambda)}\rangle^-\mathbf{1}_{\zeta})=\mathcal{T}_i(\mathbf{1}_{\zeta-\lambda}\langle{M(\lambda)}\rangle^-);
\end{equation}
\begin{eqnarray}
&&\mathcal{T}_i(\langle{M(\lambda)}\rangle^+\mathbf{1}_{\zeta})\mathcal{T}_i(\mathbf{1}_{\zeta'}\langle{M(\lambda')}\rangle^-)
=\delta_{\zeta,\zeta'}\mathcal{T}_i(\langle{M(\lambda)}\rangle^+\mathbf{1}_{\zeta}\langle{M(\lambda')}\rangle^-)\nonumber\\
&&\mathcal{T}_i(\langle{M(\lambda)}\rangle^-\mathbf{1}_{\zeta})\mathcal{T}_i(\mathbf{1}_{\zeta'}\langle{M(\lambda')}\rangle^+)
=\delta_{\zeta,\zeta'}\mathcal{T}_i(\langle{M(\lambda)}\rangle^-\mathbf{1}_{\zeta}\langle{M(\lambda')}\rangle^+).
\end{eqnarray}
\end{proposition}

For the proof of other relations, we first give some lemmas.

\begin{lemma}\label{5A}
For any $\lambda\in\mathcal{P}$ and $m\in\mathbb{N}$, we have
\begin{equation}
\mathcal{T}_i(u_i^{+(m)}\mathbf{1}_{\zeta})\mathcal{T}_i(\mathbf{1}_{\zeta'}\langle M(\lambda)\rangle^+)=\mathcal{T}_i(\delta_{\zeta,\zeta'}\mathbf{1}_{\zeta+\alpha}(u_i^{(m)}\langle M(\lambda)\rangle)^+).\label{5.1}
\end{equation}
\end{lemma}
\begin{prof}
We write $V_{\lambda}=V_{\lambda_0}\oplus tV_{i}$ as above, then
\begin{eqnarray*}
&&\mathcal{T}_i(\mathbf{1}_{\zeta+m\alpha_i}(u_i^{(m)}\langle M(\lambda)\rangle)^+)\\
&=&v^{\langle\lambda_0,ti\rangle}\mathcal{T}_i(\mathbf{1}_{\zeta+m\alpha_i}(u_i^{(m)}u_i^{(t)}\langle M(\lambda_0)\rangle)^+)\\
&=&v^{\langle\lambda_0,ti\rangle}\left[\begin{array}{c}
                                         s+t \\
                                         m
                                       \end{array}
\right]_{v_i}\mathcal{T}_i(\mathbf{1}_{\zeta+m\alpha_i}(u_i^{(m+t)}\langle M(\lambda_0)\rangle)^+)\\
&=&v^{\langle\lambda_0,ti\rangle}\left[\begin{array}{c}
                                         s+t \\
                                         m
                                       \end{array}
\right]_{v_i}
v^{-\langle\lambda_0,(m+t)i\rangle}\mathcal{T}_i(\mathbf{1}_{\zeta+m\alpha_i}(v^{\langle\lambda_0,(m+t)i\rangle}u_i^{(m+t)}\langle M(\lambda_0)\rangle)^+)\\
&=&(-1)^{m+t}\left[\begin{array}{c}
                                         s+t \\
                                         m
                                       \end{array}
\right]_{v_i}v^{r}\mathbf{1}_{s_i\zeta-m\alpha_i}u_i^{-(m+t)}\langle M(\sigma_i^+\lambda_0)\rangle^+,
\end{eqnarray*}
where
\begin{eqnarray*}
r&=&\langle\lambda_0,ti\rangle-\langle\lambda_0,(m+t)i\rangle+\langle\lambda_0,(m+t)i\rangle\\
&&(t+m)^2\varepsilon_i+t\varepsilon_i+m\varepsilon_i-(\zeta+m\alpha_i,(t+m)\alpha_i)\\
&=&\langle\lambda_0,ti\rangle+(t+m)^2\varepsilon_i+t\varepsilon_i+m\varepsilon_i-(\zeta,(t+m)\alpha_i)-2m(t+m)\varepsilon_i\\
&=&\langle\lambda_0,ti\rangle-m^2\varepsilon_i+t^2\varepsilon_i+t\varepsilon_i+m\varepsilon_i-(\zeta,(t+m)\alpha_i).
\end{eqnarray*}
While
\begin{eqnarray*}
&&\mathcal{T}_i(u_i^{+(m)}\mathbf{1}_{\zeta})\mathcal{T}_i(\mathbf{1}_{\zeta}\langle M(\lambda)\rangle^+)\\
&=&(-1)^{m}v^{r_1}\mathbf{1}_{s_i\zeta-m\alpha_i}u_i^{-(m)}(-1)^{t}v^{r_2}u_i^{-(t)}\langle M(\sigma_i^+\lambda_0)\rangle^+\\
&=&(-1)^{m+t}\left[\begin{array}{c}
                                         s+t \\
                                         m
                                       \end{array}
\right]_{v_i}v^{r_1+r_2}\mathbf{1}_{s_i\zeta-m\alpha_i}u_i^{-(m+t)}\langle M(\sigma_i^+\lambda_0)\rangle^+,
\end{eqnarray*}
where $r_1=-m^2\varepsilon_i+m\varepsilon_i-(\zeta,m\alpha_i)$ and $r_2=\langle\lambda_0,ti\rangle+t^2\varepsilon_i+t\varepsilon_i-(\zeta,t\alpha_i)$. Clearly, $r_1+r_2=r$. Hence we have formula (\ref{5.1}) in Lemma \ref{5A}.
\end{prof}

\begin{lemma}\label{5B}
For any $\lambda\in\mathcal{P}$, we have
\begin{eqnarray}
&&-(u_i^-\langle M(\lambda)\rangle^+-\langle M(\lambda)\rangle^+u_i^-)\mathbf{1}_{\zeta}\nonumber\\
&=&\frac{v_i}{a_i}(v^{(\zeta,\alpha_i)}(r_i(\langle M(\lambda)\rangle))^+-v^{(\zeta+\lambda-\alpha_i,-\alpha_i)}(r'_i(\langle M(\lambda)\rangle))^+)\mathbf{1}_{\zeta}\label{5.2}
\end{eqnarray}
and
\begin{eqnarray}
&&-(\langle M(\lambda)\rangle^-u_i^+-u_i^+\langle M(\lambda)\rangle^-)\mathbf{1}_{\zeta}\nonumber\\
&=&\frac{v_i}{a_i}(v^{(\zeta-\lambda+\alpha_i,\alpha_i)}(r'_i(\langle M(\lambda)\rangle))^--v^{(\zeta,-\alpha_i)}(r_i(\langle M(\lambda)\rangle))^-)\mathbf{1}_{\zeta}.\label{5.3}
\end{eqnarray}
\end{lemma}
\begin{prof}
Recall the relation (\ref{1})
\begin{eqnarray*}
&&\sum_{\alpha,\alpha'\in\mathcal{P}}v^{\langle\alpha',\alpha\rangle+(\alpha,\alpha)+(\zeta,-\alpha)}\frac{a_{\alpha'}}{a_{\lambda'}}g_{\alpha'\alpha}^{\lambda'}(-1)^{tr\alpha'}v^{m(\alpha')}\langle M(\alpha')\rangle^-\mathbf{1}_{\zeta+\alpha'}(r'_{\alpha}(\langle M(\lambda)\rangle))^+\nonumber\\
&&=\sum_{\alpha,\beta\in\mathcal{P}}v^{\langle\alpha,\beta\rangle+(\beta,\beta)+(\zeta,\beta)}\frac{a_{\alpha}}{a_{\lambda}}g_{\alpha\beta}^{\lambda}(-1)^{tr(\lambda'-\beta)}v^{m(\lambda'-\beta)}\langle M(\alpha)\rangle^+\mathbf{1}_{\zeta-\alpha}(r_{\beta}(\langle M(\lambda')\rangle))^-
\end{eqnarray*}
in the definition of $\dot{\mathcal{H}}^{\ast}_q(\Lambda)$.
Let $\lambda'=i$ in the above relation. We can get formula (\ref{5.2}). Similarly, let $\lambda=i$ and $\lambda'=\lambda$. We get formula (\ref{5.3}).
\end{prof}

\begin{lemma}\label{5C}
For any $\beta\in\mathcal{P}$ and $m\in\mathbb{N}$, we have
\begin{equation}
\mathcal{T}_i(\langle M(\beta)\rangle^+\mathbf{1}_{\zeta})\mathcal{T}_i(\mathbf{1}_{\zeta'}u_i^{+(m)})=\mathcal{T}_i(\delta_{\zeta,\zeta'}\mathbf{1}_{\zeta+\beta}(\langle M(\beta)\rangle u_i^{(m)})^+).\label{5.4}
\end{equation}
\end{lemma}
\begin{prof}
From the definition of $\mathcal{T}_i$, we only need to prove
\begin{displaymath}
\mathcal{T}_i(\langle M(\beta)\rangle^+\mathbf{1}_{\zeta})\mathcal{T}_i(\mathbf{1}_{\zeta}u_i^{+(m)})=\mathcal{T}_i(\mathbf{1}_{\zeta+\beta}(\langle M(\beta)\rangle u_i^{(m)})^+).
\end{displaymath}
By Lemma \ref{5A}, it suffices to prove the lemma for the case $V_{\beta}$ does not contain $V_i$ as a direct summand. So we assume that $V_i$ is not a direct summand of $V_{\beta}$.

First we have (\cite{XY})
\begin{displaymath}
\langle M(\beta)\rangle u_i=v^{(i,\beta)}u_i\langle M(\beta)\rangle+v^{-\langle{i,\beta}\rangle}\sum_{\alpha\neq\beta\oplus{i}}g_{\beta{i}}^{\alpha}\langle M(\alpha)\rangle.
\end{displaymath}
Therefore
\begin{eqnarray*}
&&\mathcal{T}_i(\mathbf{1}_{\zeta+\beta}(\langle M(\beta)\rangle u_i)^+)\\
&=&v^{(i,\beta)}\mathcal{T}_i(\mathbf{1}_{\zeta+\beta}u_i^+)\mathcal{T}_i(\langle M(\beta)\rangle^+\mathbf{1}_{\zeta-\alpha_i})\\
&&+v^{-\langle{i,\beta}\rangle}\sum_{\alpha\neq\beta\oplus{i}}g_{\beta{i}}^{\alpha}\mathcal{T}_i(\mathbf{1}_{\zeta+\beta}\langle M(\alpha)\rangle^+)\\
&=&-v^{(i,\beta)}v^{2\varepsilon_i}v^{-(\zeta+\beta,\alpha_i)}u_i^-\langle M(\sigma_i^+\beta)\rangle^+\mathbf{1}_{s_i(\zeta-\alpha_i)}\\
&&+v^{-\langle{i,\beta}\rangle}\sum_{\alpha\neq\beta\oplus{i}}g_{\beta{i}}^{\alpha}\langle M(\sigma_i^+\alpha)\rangle^+\mathbf{1}_{s_i(\zeta-\alpha_i)}.
\end{eqnarray*}
In the computation above, we use the fact that if $g_{\beta i}^{\alpha}\neq0$ and $V_{\alpha}\neq V_{\beta}\oplus V_i$, then $V_{\alpha}$ contains no direct summand isomorphic to $V_i$.
On the other hand,
\begin{eqnarray*}
&&\mathcal{T}_i(\langle M(\beta)\rangle^+\mathbf{1}_{\zeta})\mathcal{T}_i(\mathbf{1}_{\zeta}u_i^+)\\
&=&-v^{2\varepsilon_i}v^{-(\zeta,\alpha_i)}\langle M(\sigma_i^+\beta)\rangle^+\mathbf{1}_{s_i\zeta}u_i^-.
\end{eqnarray*}
Thus, to prove
\begin{displaymath}
\mathcal{T}_i(\langle M(\beta)\rangle^+\mathbf{1}_{\zeta})\mathcal{T}_i(\mathbf{1}_{\zeta}u_i^{+})=\mathcal{T}_i(\mathbf{1}_{\zeta+\beta}(\langle M(\beta)\rangle u_i)^+),
\end{displaymath}
we only need to prove
\begin{eqnarray*}
&&-v^{2\varepsilon_i}v^{-(\zeta,\alpha_i)}\langle M(\sigma_i^+\beta)\rangle^+\mathbf{1}_{s_i\zeta}u_i^-\\
&=&-v^{(i,\beta)}v^{2\varepsilon_i}v^{-(\zeta+\beta,\alpha_i)}u_i^-\langle M(\sigma_i^+\beta)\rangle^+\mathbf{1}_{s_i(\zeta-\alpha_i)}\\
&&+v^{-\langle{i,\beta}\rangle}\sum_{\alpha\neq\beta\oplus{i}}g_{\beta{i}}^{\alpha}\langle M(\sigma^+_i\alpha)\rangle^+\mathbf{1}_{s_i(\zeta-\alpha_i)}.
\end{eqnarray*}
It is sufficient to prove that
\begin{eqnarray*}
&&\langle M(\sigma_i^+\beta)\rangle^+u_i^-\mathbf{1}_{s_i(\zeta-\alpha_i)}-u_i^-\langle M(\sigma^+_i\beta)\rangle^+\mathbf{1}_{s_i(\zeta-\alpha_i)}\\
&=&-v^{-\langle{i,\beta}\rangle}v^{-2\varepsilon_i}v^{(\zeta,\alpha_i)}\sum_{\alpha\neq\beta\oplus{i}}g_{\beta{i}}^{\alpha}\langle M(\sigma_i^+\alpha)\rangle^+\mathbf{1}_{s_i(\zeta-\alpha_i)}.
\end{eqnarray*}
In rep-$\mathcal{S}$, $V_i$ is a simple injective and $V_{\sigma^+_i\beta}\in\textrm{rep-$\sigma_i\mathcal{S}$}$, so $g_{\gamma\sigma^+_i\alpha}^{\sigma^+_i\beta}=0$ for all $V_{\gamma}\in\textrm{rep-$\sigma_i\mathcal{S}$}$. By Lemma \ref{5B} we have
\begin{eqnarray*}
&&\langle M(\sigma^+_i\beta)\rangle^+u_i^-\mathbf{1}_{s_i(\zeta-\alpha_i)}-u_i^-\langle M(\sigma^+_i\beta)\rangle^+\mathbf{1}_{s_i(\zeta-\alpha_i)}\\
&=&\frac{v_i}{a_i}(v^{(s_i(\zeta-\alpha_i),\alpha_i)}(r_i(\langle M(\sigma^+_i\beta)\rangle))^+\\
&&-v^{((s_i(\zeta-\alpha_i)+s_i\beta-\alpha_i,\alpha_i)}(r'_i(\langle M(\sigma^+_i\beta)\rangle))^+)\mathbf{1}_{s_i(\zeta-\alpha_i)}\\
&=&\frac{v_i}{a_i}(v^{-(\zeta,\alpha_i)+2\varepsilon_i}(r_i(\langle M(\sigma^+_i\beta)\rangle))^+\\
&&-v^{(\zeta,\alpha_i)+(\beta,\alpha_i)}(r'_i(\langle M(\sigma^+_i\beta)\rangle))^+)\mathbf{1}_{s_i(\zeta-\alpha_i)}\\
&=&-\frac{1}{a_i}v^{(\zeta,\alpha_i)+(\beta,\alpha_i)+\varepsilon_i}(r'_i(\langle M(\sigma^+_i\beta)\rangle))^+\mathbf{1}_{s_i(\zeta-\alpha_i)}\\
&=&-\frac{1}{a_i}v^{(\zeta,\alpha_i)+(\beta,\alpha_i)+\varepsilon_i}\sum_{\alpha}\frac{a_{\sigma^+_i\alpha}a_i}{a_{\sigma^+_i\beta}}v^{\langle{i,\sigma_i^+\alpha}\rangle+(i,\sigma_i^+\alpha)}g_{i\sigma^+_i\alpha}^{\sigma^+_i\beta}\langle M(\sigma^+_i\alpha)\rangle^+\mathbf{1}_{s_i(\zeta-\alpha_i)}\\
&=&-v^{(\zeta,\alpha_i)+(\beta,\alpha_i)+\varepsilon_i}\sum_{\alpha}v^{\langle{i,\sigma_i^+\alpha}\rangle+(i,\sigma_i^+\alpha)}g_{\beta{i}}^{\alpha}\langle M(\sigma^+_i\alpha)\rangle^+\mathbf{1}_{s_i(\zeta-\alpha_i)}\\
&=&-v^{(\zeta,\alpha_i)-2\varepsilon_i-\langle{i,\beta}\rangle}\sum_{\alpha}g_{\beta{i}}^{\alpha}\langle M(\sigma^+_i\alpha)\rangle^+\mathbf{1}_{s_i(\zeta-\alpha_i)}\\
&=&-v^{(\zeta,\alpha_i)-2\varepsilon_i-\langle{i,\beta}\rangle}\sum_{\alpha}g_{\beta{i}}^{\alpha}\langle M(\sigma^+_i\alpha)\rangle^+\mathbf{1}_{s_i(\zeta-\alpha_i)}.
\end{eqnarray*}
In the computation, we use the following formula
\begin{displaymath}
g_{\beta{i}}^{\alpha}=\frac{a_{\alpha}}{a_{\beta}}g_{i\sigma^+_i\alpha}^{\sigma^+_i\beta}
\end{displaymath}
for $i\in I$ be a sink and $V_{\alpha},V_{\beta}\in\textrm{rep-$\mathcal{S}\langle i\rangle$}$.

Then by induction, we get the formula (\ref{5.4}).
\end{prof}

%

\begin{proposition}\label{proposition5A}
For $\alpha,\beta\in\mathcal{P}$, we have
\begin{equation}
\mathcal{T}_i(\langle M(\alpha)\rangle^+\mathbf{1}_{\zeta})\mathcal{T}_i(\mathbf{1}_{\zeta'}\langle M(\beta)\rangle^+)=\mathcal{T}_i(\delta_{\zeta,\zeta'}\mathbf{1}_{\zeta+\alpha}(\langle M(\alpha)\rangle\langle M(\beta)\rangle)^+).\label{5.5}
\end{equation}
\end{proposition}
\begin{prof}
By Lemma \ref{5A} and Lemma \ref{5C}, we can assume that $V_{\alpha}$ and $V_{\beta}$ do not contain $V_i$ as a direct summand.
In \cite{Ringel PBW-bases of quantum groups}, Ringel points that $\sigma_i^+$ induces an $\mathbb{Q}(v)$-algebra isomorphism from $\mathcal{H}^{\ast}_q(\Lambda)\langle{i}\rangle$ to $\mathcal{H}^{\ast}_q(\sigma_i\Lambda)\langle{i}\rangle$ mapping $\langle{M(\alpha)}\rangle$ to $\langle{M(\sigma_i^+\alpha)}\rangle$, where
$\mathcal{H}^{\ast}_q(\Lambda)\langle{i}\rangle$ is the subalgebra generated by $\langle{M(\alpha)}\rangle$ with $V_{\alpha}\in\textrm{rep-$\mathcal{S}\langle{i}\rangle$}$. Hence we prove formula (\ref{5.5}).
\end{prof}

Similarly, we have

\begin{proposition}\label{proposition5A'}
For $\alpha,\beta\in\mathcal{P}$, we have
\begin{equation}
\mathcal{T}_i(\langle M(\alpha)\rangle^-\mathbf{1}_{\zeta})\mathcal{T}_i(\mathbf{1}_{\zeta'}\langle M(\beta)\rangle^-)=\mathcal{T}_i(\delta_{\zeta,\zeta'}\mathbf{1}_{\zeta+\alpha}(\langle M(\alpha)\rangle\langle M(\beta)\rangle)^-).\label{5.6}
\end{equation}
\end{proposition}

Then the most difficult defining relation (\ref{2}) should be verified, that is, for an element $y\in\dot{\mathcal{H}}_q^{\ast}(\Lambda)$, which can be writen as $$y=\sum_{x,x',\zeta} x^+\mathbf{1}_{\zeta}x'^-$$ and $$y=\sum_{x,x',\zeta} x^-\mathbf{1}_{\zeta}x'^+,$$ we should verify that
$$\sum_{x,x',\zeta} T_{i}(x^+\mathbf{1}_{\zeta}x'^-)=\sum_{x,x',\zeta} T_{i}(x^-\mathbf{1}_{\zeta}x'^+).$$

\begin{proposition}\label{proposition5C}
For any $\lambda,\lambda'\in\mathcal{P}$, we have
\begin{eqnarray}
\sum_{\alpha,\alpha'\in\mathcal{P}}v^{\langle\alpha',\alpha\rangle+(\alpha,\alpha)+(\zeta,-\alpha)}\frac{a_{\alpha'}}{a_{\lambda'}}g_{\alpha'\alpha}^{\lambda'}(-1)^{tr\alpha'}v^{m(\alpha')}
\mathcal{T}_i\left(\langle M(\alpha')\rangle^-\mathbf{1}_{\zeta+\alpha'}(r'_{\alpha}(\langle M(\lambda)\rangle))^+\right)=\nonumber\\
\sum_{\alpha,\beta\in\mathcal{P}}v^{\langle\alpha,\beta\rangle+(\beta,\beta)+(\zeta,\beta)}\frac{a_{\alpha}}{a_{\lambda}}g_{\alpha\beta}^{\lambda}(-1)^{tr(\lambda'-\beta)}v^{m(\lambda'-\beta)}
\mathcal{T}_i\left(\langle M(\alpha)\rangle^+\mathbf{1}_{\zeta-\alpha}(r_{\beta}(\langle M(\lambda')\rangle))^-\right).\label{5.7}
\end{eqnarray}
\end{proposition}

\begin{prof}
By Proposition \ref{proposition5A} and Proposition \ref{proposition5A'}, we may assume that $V_{\lambda}$ and $V_{\lambda'}$ contain no direct summand isomorphic to $V_i$. Then $V_{\alpha}$ and $V_{\alpha'}$ also contain no direct summand isomorphic to $V_i$.

Let
\begin{displaymath}
L=\sum_{\alpha,\alpha'\in\mathcal{P}}v^{\langle\alpha',\alpha\rangle+(\alpha,\alpha)+(\zeta,-\alpha)}\frac{a_{\alpha'}}{a_{\lambda'}}g_{\alpha'\alpha}^{\lambda'}(-1)^{tr\alpha'}v^{m(\alpha')}\langle M(\alpha')\rangle^-\mathbf{1}_{\zeta+\alpha'}(r'_{\alpha}(\langle M(\lambda)\rangle))^+;
\end{displaymath}
and
\begin{displaymath}
R=\sum_{\alpha,\beta\in\mathcal{P}}v^{\langle\alpha,\beta\rangle+(\beta,\beta)+(\zeta,\beta)}\frac{a_{\alpha}}{a_{\lambda}}g_{\alpha\beta}^{\lambda}(-1)^{tr(\lambda'-\beta)}v^{m(\lambda'-\beta)}\langle M(\alpha)\rangle^+\mathbf{1}_{\zeta-\alpha}(r_{\beta}(\langle M(\lambda')\rangle))^-.
\end{displaymath}

First consider $L$. We have
\begin{eqnarray*}
L&=&\mathbf{1}_{\zeta}\sum_{\alpha,\alpha',\beta\in\mathcal{P}}v^{\langle{\lambda',\alpha}\rangle+(\zeta,-\alpha)+\langle{\alpha,\lambda}\rangle+(\alpha,\beta)}
\frac{a_{\alpha'}a_{\alpha}a_{\beta}}{a_{\lambda'}a_{\lambda}}g_{\alpha'\alpha}^{\lambda'}g^{\lambda}_{\alpha\beta}(-1)^{tr\alpha'}v^{m(\alpha')}
\langle{M(\alpha')}\rangle^-\langle M(\beta)\rangle)^+\\
&=&\mathbf{1}_{\zeta}\sum_{\alpha,\alpha',\beta\in\mathcal{P}}A_1B_1\langle{M(\alpha')}\rangle^-\langle M(\beta)\rangle)^+\\
\end{eqnarray*}
where $A_1=v^{\langle{\lambda',\alpha}\rangle+(\zeta,-\alpha)+\langle{\alpha,\lambda}\rangle+(\alpha,\beta)}(-1)^{tr\alpha'}v^{m(\alpha')}$ and $B_1=\frac{a_{\alpha'}a_{\alpha}a_{\beta}}{a_{\lambda'}a_{\lambda}}g_{\alpha'\alpha}^{\lambda'}g^{\lambda}_{\alpha\beta}$.

Now assume $V_{\beta}=V_{\beta'}\oplus{tV_i}$, where $V_{\beta'}$ contains no direct summand isomorphic to $V_i$. Then we have $\langle{M(\beta)}\rangle=v^{\langle{\beta',ti}\rangle}u_i^{(t)}\langle{M(\beta')}\rangle$.

Then
\begin{eqnarray*}
&&\mathcal{T}_i(L)\\
&=&\mathbf{1}_{s_i\zeta}\sum_{\alpha,\alpha',\beta\in\mathcal{P}}A_1B_1
\mathcal{T}_i(\langle{M(\alpha')}\rangle^-\mathbf{1}_{\zeta+\alpha'}\langle M(\beta)\rangle)^+)\\
&=&\mathbf{1}_{s_i\zeta}\sum_{\alpha,\alpha',\beta'\in\mathcal{P},t}A_1B_1
(-1)^{t-\alpha'(h_i)}v^{t^2\varepsilon_i+t\varepsilon_i+\langle{\beta',ti}\rangle-(\zeta+\alpha',t\alpha_i)-(\alpha',i)}\\
&&\langle{M(\sigma^+_i \alpha')}\rangle^-u_i^{-(t)}\langle{M(\sigma^+_i\beta')}\rangle)^+\\
&=&\mathbf{1}_{s_i\zeta}\sum_{\alpha,\alpha',\beta'\in\mathcal{P},t}A_1B_1A_2\langle{M(\sigma^+_i \alpha')}\rangle^-u_i^{-(t)}\langle{M(\sigma^+_i\beta')}\rangle)^+\\
\end{eqnarray*}
where $A_2=(-1)^{t-\alpha'(h_i)}v^{t^2\varepsilon_i+t\varepsilon_i+\langle{\beta',ti}\rangle-(\zeta+\alpha',t\alpha_i)-(\alpha',i)}$.

Since $i$ is a source of $\sigma_iQ$ and $V_{\alpha'}$ contains no direct summand isomorphic to $V_i$, $\langle{M(\sigma^+_i\alpha'\oplus{ti})}\rangle=v^{\langle{ti,\alpha'}\rangle}\langle{M(\sigma^+_i\alpha')}\rangle u_i^{(t)}$.

Hence we have
\begin{displaymath}
\mathcal{T}_i(L)=\mathbf{1}_{s_i\zeta}\sum_{\alpha,\alpha',\beta'\in\mathcal{P},t}A_1B_1A_2A_3
\langle{M(\sigma^+_i\alpha'\oplus{ti})}\rangle^-\langle{M(\sigma^+_i\beta')}\rangle)^+
\end{displaymath}
where $A_3=v^{-\langle{ti,\alpha'}\rangle}$.

Then we compute $B_1$.

If $i$ is a sink and $V_{\alpha},V_{\beta}$ contain no direct summand isomorphic to $V_i$, then $g_{\alpha,\beta\oplus{ti}}^{\lambda}=\sum_{\gamma}g^{\gamma}_{\alpha{ti}}g^{\lambda}_{\gamma\beta}$. If $i$ is a source and $V_{\alpha},V_{\beta}$ contain no direct summand isomorphic to $V_i$, then $g_{\alpha\oplus{ti},\beta}^{\lambda}=\sum_{\gamma}g^{\gamma}_{{ti}\beta}g^{\lambda}_{\alpha\gamma}$.

Since $V_{\alpha}$ and $V_{\beta'}$ contain no direct summand isomorphic to $V_i$, we have
\begin{displaymath}
g_{\alpha\beta}^{\lambda}=\sum_{\gamma}g^{\gamma}_{\alpha{ti}}g^{\lambda}_{\gamma\beta'}.
\end{displaymath}
Note that (\cite{XY})
\begin{displaymath}
a_{\beta}=v^{2\langle{ti,\beta'}\rangle}a_{\beta'}a_{ti}, a_{\sigma_i^+\alpha'\oplus{ti}}=v^{2\langle{ti,\alpha'}\rangle}a_{\alpha'}a_{ti}.
\end{displaymath}
Then
\begin{eqnarray*}
B_1&=&\frac{a_{\alpha'}a_{\alpha}a_{\beta}}{a_{\lambda'}a_{\lambda}}g_{\alpha'\alpha}^{\lambda'}g^{\lambda}_{\alpha\beta}\\
&=&\sum_{\gamma}v^{2\langle{ti,\beta'}\rangle}\frac{a_{\alpha'}a_{\alpha}a_{\beta'}a_{ti}}{a_{\lambda'}a_{\lambda}}
g_{\alpha'\alpha}^{\lambda'}g^{\gamma}_{\alpha{ti}}g^{\lambda}_{\gamma\beta'}.\\
\end{eqnarray*}
We may assume $V_{\gamma}$ contains no direct summand isomorphic to $V_i$. Hence we have
\begin{displaymath}
a_{\alpha}g_{\gamma}^{\alpha{ti}}=a_{\gamma}g^{\sigma^+_i\alpha}_{ti\sigma^+_i\gamma}.
\end{displaymath}
Then
\begin{eqnarray*}
B_1
&=&\sum_{\gamma}v^{2\langle{ti,\beta'}\rangle}\frac{a_{\alpha'}a_{\gamma}a_{\beta'}a_{ti}}{a_{\lambda'}a_{\lambda}}
g_{\alpha'\alpha}^{\lambda'}g^{\sigma^+_i\alpha}_{ti\sigma^+_i\gamma}g^{\lambda}_{\gamma\beta'}\\
&=&\sum_{\gamma}v^{2\langle{ti,\beta'}\rangle}\frac{a_{\sigma^+_i\alpha'}a_{\sigma^+_i\gamma}a_{\sigma^+_i\beta'}a_{ti}}{a_{\sigma^+_i\lambda'}a_{\sigma^+_i\lambda}}
g_{\sigma^+_i\alpha'\sigma^+_i\alpha}^{\sigma^+_i\lambda'}g^{\sigma^+_i\alpha}_{ti\sigma^+_i\gamma}g^{\sigma^+_i\lambda}_{\sigma^+_i\gamma\sigma^+_i\beta'}\\
&=&\sum_{\gamma}v^{2\langle{ti,\beta'}\rangle}v^{2\langle{\alpha',ti}\rangle}\frac{a_{\sigma^+_i\alpha'\oplus{ti}}a_{\sigma^+_i\gamma}a_{\sigma^+_i\beta'}}{a_{\sigma^+_i\lambda'}a_{\sigma^+_i\lambda}}
g_{\sigma^+_i\alpha'\sigma^+_i\alpha}^{\sigma^+_i\lambda'}g^{\sigma^+_i\alpha}_{ti\sigma^+_i\gamma}g^{\sigma^+_i\lambda}_{\sigma^+_i\gamma\sigma^+_i\beta'}\\
&=&\sum_{\gamma}A_4\frac{a_{\sigma^+_i\alpha'\oplus{ti}}a_{\sigma^+_i\gamma}a_{\sigma^+_i\beta'}}{a_{\sigma^+_i\lambda'}a_{\sigma^+_i\lambda}}
g_{\sigma^+_i\alpha'\sigma^+_i\alpha}^{\sigma^+_i\lambda'}g^{\sigma^+_i\alpha}_{ti\sigma^+_i\gamma}g^{\sigma^+_i\lambda}_{\sigma^+_i\gamma\sigma^+_i\beta'}\\
&=&\sum_{\gamma}A_4\frac{a_{\sigma^+_i\alpha'\oplus{ti}}a_{\sigma^+_i\gamma}a_{\sigma^+_i\beta'}}{a_{\sigma^+_i\lambda'}a_{\sigma^+_i\lambda}}
g_{\sigma^+_i\alpha'\oplus{ti},\sigma^+_i\gamma}^{\sigma^+_i\lambda'}g^{\sigma^+_i\lambda}_{\sigma^+_i\gamma\sigma^+_i\beta'}\\
\end{eqnarray*}
where $A_4=v^{2\langle{ti,\beta'}\rangle}v^{2\langle{\alpha',ti}\rangle}$.

Then we compute $A=A_1A_2A_3A_4$.
\begin{eqnarray*}
A&=&A_1A_2A_3A_4\\
&=&v^{\langle{\lambda',\alpha}\rangle+(\zeta,-\alpha)+\langle{\alpha,\lambda}\rangle+(\alpha,\beta)}(-1)^{tr\alpha'}v^{m(\alpha')}\\
&&(-1)^{t-\alpha'(h_i)}v^{t^2\varepsilon_i+t\varepsilon_i+\langle{\beta',ti}\rangle-(\zeta+\alpha',t\alpha_i)-(\alpha',i)}\\
&&v^{-\langle{ti,\alpha'}\rangle}\\
&&v^{2\langle{ti,\beta'}\rangle}v^{2\langle{\alpha',ti}\rangle}\\
&=&(-1)^{tr(\sigma_i^+(\alpha'))+t}
v^{(\zeta,-\alpha-t\alpha_i)+\langle{\sigma_i^+(\lambda'),\sigma_i^+(\gamma)}\rangle+\langle{\sigma_i^+(\gamma),\sigma_i^+(\lambda)}\rangle+(\sigma_i^+(\gamma),\sigma_i^+(\beta'))+m(\sigma_i^+(\alpha'))+t\varepsilon_i}.
\end{eqnarray*}
Let $\mu_1=\sigma^+_i\gamma$, $\mu_2=\sigma^+_i\beta'$ and $\mu_3=\sigma^+_i\alpha'\oplus{ti}$.
Hence we have
\begin{eqnarray*}
\mathcal{T}_i(L)&=&\mathbf{1}_{s_i\zeta}\sum_{\mu_1,\mu_2,\mu_3\in\mathcal{P}}(-1)^{tr\mu_3}v^{m(\mu_3)}
v^{(s_i\zeta,-\mu_1)+\langle{\sigma_i^+(\lambda'),\mu_1}\rangle+\langle{\mu_1,\sigma_i^+(\lambda)}\rangle+(\mu_1,\mu_2)}\\
&&\frac{a_{\mu_3}a_{\mu_1}a_{\mu_2}}{a_{\sigma^+_i\lambda'}a_{\sigma^+_i\lambda}}
g_{\mu_3,\mu_1}^{\sigma^+_i\lambda'}g^{\sigma^+_i\lambda}_{\mu_1\mu_2}
\langle{M(\mu_3)}\rangle^-\langle{M(\mu_2)}\rangle)^+\\
&=&\mathbf{1}_{s_i\zeta}\sum_{\mu_1,\mu_2,\mu_3\in\mathcal{P}}(-1)^{tr\mu_3}v^{m(\mu_3)}
v^{(s_i\zeta,-\mu_1)+\langle{\mu_1+\mu_3,\mu_1}\rangle+\langle{\mu_1,\mu_1+\mu_2}\rangle+(\mu_1,\mu_2)}\\
&&\frac{a_{\mu_3}a_{\mu_1}a_{\mu_2}}{a_{\sigma^+_i\lambda'}a_{\sigma^+_i\lambda}}
g_{\mu_3,\mu_1}^{\sigma^+_i\lambda'}g^{\sigma^+_i\lambda}_{\mu_1\mu_2}
\langle{M(\mu_3)}\rangle^-\langle{M(\mu_2)}\rangle)^+\\
&=&\mathbf{1}_{s_i\zeta}\sum_{\mu_1,\mu_3\in\mathcal{P}}(-1)^{tr\mu_3}v^{m(\mu_3)}
v^{(s_i\zeta,-\mu_1)+\langle{\mu_3,\mu_1}\rangle+(\mu_1,\mu_1)}\\
&&\frac{a_{\mu_3}}{a_{\sigma^+_i\lambda'}}g_{\mu_3,\mu_1}^{\sigma^+_i\lambda'}
\langle{M(\mu_3)}\rangle^-(r'_{\mu_1}(\langle{M(\sigma^+_i\lambda)}\rangle))^+.
\end{eqnarray*}

Similarly we have

\begin{eqnarray*}
\mathcal{T}_i(R)
&=&\mathbf{1}_{s_i\zeta}\sum_{\mu_4,\mu_5\in\mathcal{P}}(-1)^{tr(s_i\lambda'-\mu_5)}v^{m(s_i\lambda'-\mu_5)}
v^{(s_i\zeta,\mu_5)+\langle{\mu_4,\mu_5}\rangle+(\mu_5,\mu_5)}\\
&&\frac{a_{\mu_4}}{a_{\sigma^+_i\lambda}}
g^{\sigma^+_i\lambda}_{\mu_4\mu_5}
\langle{M(\mu_4)}\rangle^+(r_{\mu_5}(\langle{M(\sigma^+_i\lambda')}\rangle))^-.
\end{eqnarray*}

By the first relation (\ref{1}) in the definition of $\dot{\mathcal{H}}_q^{\ast}(\Lambda)$, we have $\mathcal{T}_i(L)=\mathcal{T}_i(R)$.

\end{prof}

Then Proposition \ref{proposition5B}, \ref{proposition5A}, \ref{proposition5A'} and \ref{proposition5C} imply Theorem \ref{Theorem 5D}.


\begin{thebibliography}{99}

%

\bibitem{Lusztig quantum deformation}
Lusztig, G.: Quantum deformations of certain simple modules over enveloping algebras. \emph{Adv. Math.}, \textbf{70}, 237--249 (1988)

\bibitem{Lusztig quantum groups at roots of 1}
Lusztig, G.: Quantum groups at roots of $1$. \emph{Geom. Dedicata}, \textbf{35}, 89--114 (1990)

\bibitem{Lusztig introduction to quantum groups} 
Lusztig, G.: Introduction to quantum groups, Birkhauser, Boston, 1993


\bibitem{Green}
Green, J. A.: Hall algebras, hereditary algebras and quantum groups. \emph{Invent. Math.}, \textbf{120}, 361--377 (1995)

\bibitem{Ringel Hall algebras and quantum groups}
Ringel, C. M.: Hall algebras and quantum groups. \emph{Invent. Math.}, \textbf{101}, 583--592 (1990)

\bibitem{X Drinfed Double and Ringel-Green Theorem and Hall Algebras}
Xiao, J.: Drinfeld double and Ringel-Green theory of Hall algebras. \emph{J. Algebra}, \textbf{190}, 100--144 (1997)


\bibitem{SV On the double of the Hall algebra of a quiver}
Sevenhant, B., Van den Bergh, M.: On the double of the Hall algebra of a quiver. \emph{J. Algebra}, \textbf{221}, 135--160 (1999)



\bibitem{XY}
Xiao, J., Yang, S.: BGP-Reflection Functors and Lusztig's Symmetries: A Ringel-Hall Algebra Approach to Quantum Groups. \emph{J. Algebra}, \textbf{241}, 204--246 (2001)

\bibitem{DX On double Ringel¨CHall algebras}
Deng, B., Xiao, J.: On double Ringel-Hall algebras. \emph{J. Algebra}, \textbf{251}, 110--149 (2002)

\bibitem{DX Ringel¨CHall algebras and Lusztig's symmetries}
Deng, B., Xiao, J.: Ringel-Hall algebras and Lusztig's symmetries. \emph{J. Algebra}, \textbf{255}, 357--372 (2002)

\bibitem{Kac} 
Kac, V. G.: Infinite dimensional Lie algebras, Third edition, Cambridge Univ. Press, Cambridge, 1990

\bibitem{Dlab Ringel}
Dlab, V., Ringel, C. M.: Indecomposable representation of graphs and algebras. \emph{Mem. Amer. Math. Sco.}, \textbf{173}, (1976)

\bibitem{BGP}
Bernstein, I. N., Gelfand, I. M., Ponomarev, A.: Coxter functors and Gabriel's Theorem. \emph{Lispehi Math. Nauk}, \textbf{28}, 19--33 (1972)

\bibitem{SV2}
Sevenhant, B., Van den Bergh, M.: A relation between a conjecture of Kac and the structure of the Hall algebra. \emph{J. Pure Appl, Algebra}, \textbf{160}, 319--332 (2001)

%
%
%

\bibitem{Ringel PBW-bases of quantum groups}
Ringel, C. M.: PBW-bases of quantum groups. \emph{J. reine angew. Math.}, \textbf{470}, 51--88 (1996)





















\end{thebibliography}
\end{document}